\documentclass[12pt,a4paper]{amsart}
\usepackage{mathtools}
\usepackage{xcolor}
\usepackage{amsmath,amstext,amssymb,amsthm,amsfonts}
\usepackage{graphicx}
\usepackage[all]{xy} \xyoption{poly}
\usepackage{mathrsfs}  

\newcommand{\nc}{\newcommand}

\nc{\one}{\mbox{\bf 1}}
\nc{\invtensor}{\underset{\leftarrow}{\otimes}}
\nc{\const}{\operatorname{const}}

\nc{\ad}{\operatorname{ad}}

\nc{\tr}{\operatorname{tr}}

\nc{\Gr}{\mathscr{K}}
\nc{\rGr}{\operatorname{rGr}}
\nc{\atyp}{\operatorname{atyp}}
\nc{\tp}{\operatorname{top}}
\nc{\rank}{\operatorname{rank}}
\nc{\corank}{\operatorname{corank}}
\nc{\codim}{\operatorname{codim}}
\nc{\sdim}{\operatorname{sdim}}
\nc{\mult}{\operatorname{mult}}
\nc{\ds}{\operatorname{ds}}
\nc{\ext}{\operatorname{ext}}
\nc{\tail}{\operatorname{tail}}
\nc{\howl}{\operatorname{howl}}
\nc{\spn}{\operatorname{span}}
\nc{\Sym}{\operatorname{Sym}}
\nc{\core}{\operatorname{core}}
\nc{\id}{\operatorname{id}}
\nc{\Id}{\operatorname{Id}}
\nc{\Ree}{\operatorname{Re}}
\nc{\hi}{\operatorname{hi}}
\nc{\htt}{\operatorname{ht}}
\nc{\at}{\operatorname{at}}
\nc{\str}{\operatorname{str}}
\nc{\Iso}{\operatorname{Iso}}
\nc{\Ker}{\operatorname{Ker}}
\nc{\rker}{\operatorname{rKer}}
\nc{\im}{\operatorname{Im}}
\nc{\osp}{\mathfrak{osp}}
\nc{\sgn}{\operatorname{sgn}}
\nc{\F}{\operatorname{F}}
\nc{\Mod}{\operatorname{Mod}}
\nc{\DS}{\operatorname{DS}}
\nc{\Soc}{\operatorname{Soc}}
\nc{\Inj}{\operatorname{Inj}}
\nc{\Hom}{\operatorname{Hom}}
\nc{\End}{\operatorname{End}}
\nc{\supp}{\operatorname{supp}}
\nc{\Card}{\operatorname{Card}}
\nc{\Ann}{\operatorname{Ann}}
\nc{\Arc}{\operatorname{Arc}}
\nc{\arc}{\operatorname{arc}}
\nc{\Ind}{\operatorname{Ind}}
\nc{\Coind}{\operatorname{Coind}}
\nc{\wt}{\operatorname{wt}}
\nc{\hwt}{\operatorname{wt}}
\nc{\rk}{\operatorname{rank}}
\nc{\ch}{\operatorname{ch}}
\nc{\sch}{\operatorname{sch}}
\nc{\mdim}{\operatorname{mdim}}
\nc{\Stab}{\operatorname{Stab}}
\nc{\Ima}{\operatorname{Im}}
\nc{\Irr}{\operatorname{Irr}}
\nc{\Spec}{\operatorname{Spec}}
\nc{\Res}{\operatorname{Res}}
\nc{\res}{\operatorname{res}}
\nc{\Aut}{\operatorname{Aut}}
\nc{\Ext}{\operatorname{Ext}}
\nc{\Prec}{\operatorname{Prec}}
\nc{\Fract}{\operatorname{Fract}}
\nc{\gr}{\operatorname{gr}}
\nc{\diag}{\operatorname{diag}}
\nc{\deff}{\operatorname{def}}
\nc{\depth}{\operatorname{depth}}
\nc{\HC}{\operatorname{HC}}
\nc{\dpth}{\operatorname{dpth}}
\nc{\sw}{\operatorname{sw}}
\nc{\red}{\operatorname{red}}
\nc{\pari}{\operatorname{par}}
\nc{\pos}{\operatorname{pos}}

\nc{\Cl}{\mathcal{C}\ell}

\nc{\wdchi}{\widetilde{\chi}}
\nc{\wdH}{\widetilde{H}}
\nc{\wdN}{\widetilde{N}}
\nc{\wdM}{\widetilde{M}}
\nc{\wdO}{\widetilde{O}}
\nc{\wdR}{\widetilde{R}}

\nc{\wdV}{\widetilde{V}}

\nc{\wdC}{\widetilde{C}}

\nc{\zero}{\operatorname{zero}}
\nc{\nonzero}{\operatorname{nonzero}}

\nc{\Obj}{\operatorname{Obj}}
\nc{\Dglie}{\operatorname{{\mathcal D}glie}}
\nc{\Fin}{\operatorname{{\mathcal F}in}}
\nc{\pr}{\operatorname{pr}}
\nc{\Adm}{\operatorname{\mathcal{A}dm}}

\nc{\Sg}{{\cS(\fg)}}
\nc{\Shg}{{\cS(\fhg)}}
\nc{\Ug}{{\cU(\fg)}}
\nc{\Uhg}{{\cU(\fhg)}}
\nc{\Sh}{{\cS(\fh)}}
\nc{\Uh}{{\cU(\fh)}}
\nc{\Uhh}{{\cU(\fhh)}}
\nc{\Zg}{{{\mathcal{Z}}(\fg)}}

\nc{\Vir}{{\mathcal{V}ir}}
\nc{\NS}{{\mathcal{N}S}}

\nc{\tZg}{{\widetilde{\mathcal Z}({\mathfrak g})}}
\nc{\Zk}{{\mathcal Z}({\mathfrak k})}

\newcommand{\shr}{shr}

\nc{\Up}{{\mathcal U}({\mathfrak p})}
\nc{\Ah}{{\mathcal A}({\mathfrak h})}
\nc{\Ag}{{\mathcal A}({\mathfrak g})}
\nc{\Ap}{{\mathcal A}({\mathfrak p})}
\nc{\Zp}{{\mathcal Z}({\mathfrak p})}
\nc{\cR}{\mathcal R}
\nc{\cS}{\mathcal S}
\nc{\cP}{\mathcal P}
\nc{\cT}{\mathcal{T}}
\nc{\CC}{\mathcal C}
\nc{\cA}{\mathcal A}
\nc{\cU}{\mathcal U}
\nc{\cZ}{\mathcal Z}
\nc{\cM}{\mathcal M}
\nc{\cL}{\mathcal L}
\nc{\cF}{\mathcal F}
\nc{\fg}{\mathfrak g}
\nc{\cB}{\mathcal{B}}

\nc{\fo}{\mathfrak o}

\nc{\CO}{\mathcal O}
\nc{\CR}{\mathcal R}

\nc{\cK}{\mathcal{K}}
\nc{\cW}{\mathcal{W}}
\nc{\bM}{\mathbf{M}}
\nc{\bL}{\mathbf{L}}
\nc{\bN}{\mathbf{N}}

\nc{\zq}{\mathpzc q}

\nc{\fl}{\mathfrak l}
\nc{\fn}{\mathfrak n}
\nc{\fm}{\mathfrak m}
\nc{\fp}{\mathfrak p}
\nc{\fh}{\mathfrak h}
\nc{\ft}{\mathfrak t}
\nc{\fk}{\mathfrak k}
\nc{\fb}{\mathfrak b}
\nc{\fs}{\mathfrak s}
\nc{\fB}{\mathfrak B}

\nc{\vareps}{\varepsilon}
\nc{\varesp}{\varepsilon}
\nc{\veps}{\varepsilon}

\nc{\fsl}{\mathfrak{sl}}
\nc{\fgl}{\mathfrak{gl}}
\nc{\fso}{\mathfrak{so}}
\nc{\fosp}{\mathfrak{osp}}
\nc{\fsp}{\mathfrak{sp}}
\nc{\fq}{\mathfrak q}
\nc{\fsq}{\mathfrak{sq}}
\nc{\fpsq}{\mathfrak{psq}}
\nc{\fpq}{\mathfrak{pq}}

\nc{\fhg}{\hat{\fg}}
\nc{\fhn}{\hat{\fn}}
\nc{\fhh}{\hat{\fh}}
\nc{\fhb}{\hat{\fb}}
\nc{\hrho}{\hat{\rho}}

\nc{\hsl}{\hat{\fsl}}
\nc{\fpo}{\mathfrak{po}}
\nc{\dirlim}{\underset{\rightarrow}{\lim}\,}
\nc{\nen}{\newenvironment}
\nc{\ol}{\overline}
\nc{\ul}{\underline}
\nc{\ra}{\rightarrow}
\nc{\lra}{\longrightarrow}
\nc{\Lra}{\Longrightarrow}
\nc{\bo}{\bar{1}}
\nc{\Lla}{\Longleftarrow}

\nc{\Llra}{\Longleftrightarrow}

\nc{\thla}{\twoheadleftarrow}

\nc{\lang}{(}
\nc{\rang}{)}

\nc{\hra}{\hookrightarrow}

\nc{\iso}{\overset{\sim}{\lra}}

\nc{\ssubset}{\underset{\not=}{\subset}}

\nc{\vac}{|0\rangle}

\nc{\simka}{{\ \scriptscriptstyle _{{\sim}}^\text{\tiny{k}}\ }}

\nc{\Thm}[1]{Theorem~\ref{#1}}
\nc{\Prop}[1]{Proposition~\ref{#1}}
\nc{\Lem}[1]{Lemma~\ref{#1}}
\nc{\Cor}[1]{Corollary~\ref{#1}}
\nc{\Conj}[1]{Conjecture~\ref{#1}}
\nc{\Claim}[1]{Claim~\ref{#1}}
\nc{\Defn}[1]{Definition~\ref{#1}}
\nc{\Exa}[1]{Example~\ref{#1}}
\nc{\Rem}[1]{Remark~\ref{#1}}
\nc{\Note}[1]{Note~\ref{#1}}
\nc{\Quest}[1]{Question~\ref{#1}}
\nc{\Hyp}[1]{Hypoth\`ese~\ref{#1}}
\nen{thm}[1]{\label{#1}{\bf Theorem.\ } \em}{}
\nen{prop}[1]{\label{#1}{\bf Proposition.\ } \em}{}
\nen{lem}[1]{\label{#1}{\bf Lemma.\ } \em}{}
\nen{cor}[1]{\label{#1}{\bf Corollary.\ } \em}{}
\nen{conj}[1]{\label{#1}{\bf Conjecture.\ } \em}{}

\nen{claim}[1]{\label{#1}{\bf Claim.\ } \em}{}

\nen{defn}[1]{\label{#1}{\bf Definition.\ } }{}
\nen{exa}[1]{\label{#1}{\bf Example.\ } }{}
\nen{rem}[1]{\label{#1}{\em Remark.\ } }{}
\nen{note}[1]{\label{#1}{\em Note.\ } }{}
\nen{exer}[1]{\label{#1}{\em Exercise.\ } }{}
\nen{sket}[1]{\label{#1}{\em Sketch of proof.\ } }{}
\nen{quest}[1]{\label{#1}{\bf Question.\ } \em}{}

\nen{hyp}[1]{\label{#1}{\bf Hypoth\`ese.\ } \em}{}
\usepackage{tikz}

\begin{document}
\setcounter{section}{0}
\setcounter{tocdepth}{1}

\title[$\DS$-functor for the queer Lie superalgebra]{On the Duflo-Serganova functor for the queer Lie superalgebra}
\author{ M. Gorelik, A. Sherman}

\address{Weizmann Institute of Science, Rehovot, Israel}
\email{maria.gorelik@gmail.com}
\address{Ben Gurion University of the Negev, Be'er Sheva, Israel}
\email{xandersherm@gmail.com}

\date{}

\begin{abstract} We study the Duflo-Serganova functor $\DS_x$ for the queer Lie superalgebra $\fq_n$ and for all odd $x$ with $[x,x]$ semisimple.  For the case when the rank of $x$ is $1$ we give a formula for multiplicities in terms of the arch diagram attached to $\lambda$. Further, we prove that $\DS_x(L)$ is semisimple
if $L$ is a simple finite-dimensional module and $x$ is of rank $1$ satisfying $x^2=0$.
\end{abstract}


\medskip


\maketitle

\section{Introduction}
For a finite dimensional complex Lie superalgebra $\mathfrak{g}$ and an odd element $x$ satisfying $[x,x]=0$, M.~Duflo and V.~Serganova defined a functor $\DS_x: Rep(\mathfrak{g}) \to Rep(\mathfrak{g}_x)$ where $\mathfrak{g}_x:= \Ker\ad(x) / \Ima\ad(x)$. For $\fg=\mathfrak{gl}(m|n),\osp(m|2n), \fp_n, \fq_n$ the algebra $\fg_x$ is isomorphic to
$\mathfrak{gl}(m-s|n-s),\osp(m-2s|2n-2s), \fp_{n-s}, \fq_{n-2s}$
respectively where $s$ is a non-negative number called  {\em rank of $x$}.

\subsection{Previous results}\label{KMp}
Let $\fg$ be a finite-dimensional Kac-Moody superalgebra  and
let $L:=L_{\fg}(\lambda)$ be a simple finite-dimensional $\fg$-module of the highest weight $\lambda$. For a fixed $x$ of rank $r$
we denote $\DS_x$ by $\DS_r$ (we will only use this notation for Kac-Moody superalgebras). The following properties of the $\DS$-functor were obtained in~\cite{HW},\cite{GH},\cite{Gdex}:
\begin{itemize}
\item [$(a)$]
the multiplicities of irreducible constituents of $\DS_1(L)$ 
are at most $2$;

\item[$(b)$] the following are equivalent
\begin{itemize}
\item $L$ is typical;

\item $\DS_x(L)=0$ for all non-zero $x$;

\item 
$\DS_1(L)=0$;
\end{itemize}

\item[$(c)$]  $\DS_x(L)$ does not have subquotients differing
by the parity change;

\item[$(d)$] $\DS_x(L)$ is semisimple;

\item[$(e)$]
$\DS_s(L)\cong\DS_1(\DS_1(\ldots(\DS_1(L)\ldots))$.
\end{itemize}

The $\fp_n$-case was studied in~\cite{EAS}: in this case
$(a)$--$(c)$ hold and $(d)$, $(e)$ do not hold. 

\subsubsection{}
Throughout this article, $\Pi$ stands for the parity change functor.  We will use
the notation $[N:L]$
for the ``non-graded multiplicity'' of a simple module
$L$ in a finite length module $N$: this is equal to the usual 
multiplicity  if $L\cong \Pi L$, and to the sum
of multiplicites of $L$ and of $\Pi L$ if $L\not\cong \Pi L$.  

In the ``non-exceptional cases''
$\fgl(m|n),\osp(m|n)$ and $\fp_n$ the non-graded multiplicity
$[\DS_1(L_{\fg}(\lambda)): L_{\fg_x}(\nu)]$  is given in terms of
so-called ``arc  diagrams'': this multiplicity is non-zero if and only if
the arc diagram $\Arc(\nu)$ can be 
 obtained from the arc diagram  $\Arc(\lambda)$  by removing
a maximal arc; this multiplicity is $1$ for  $\fgl(m|n)$, $\fp_n$
and  is $1$ or $2$  for $\osp(m|n)$.

\subsection{Extension of the Duflo-Serganova functor}\label{gss}

In this paper  we study $\DS_x(L)$ for the case when $L=L_{\fq_n}(\lambda)$ is a simple finite-dimensional  $\fq_n$-module.  For $\fq_n$ it is important to expand our use of the $\DS$-functor to odd elements $x$ for which $[x,x]$ is a semisimple element of $\fg_{\ol{0}}$; this is especially true because of the Cartan subalgebra admitting a nontrivial odd part.  

Write $x^2:=\frac{1}{2}[x,x]\in\fg_{\ol{0}}$.   Recall that if $x\in\fg_{\ol{1}}$ has that $x^2$ acts semisimply on a $\fg$-module $M$, we may define the Lie superalgebra
\[
M_x:=\frac{\ker(x:M^{x^2})}{\im(x:M^{x^2})},
\]
where $M^{x^2}$ is the subspace of $M$ killed by $x^2$.  In particular if $\ad(x^2)$ acts semisimply, we obtain a Lie superalgebra $\fg_x$, and $M_x$ will naturally be a module over $\fg_x$.  Thus we will be interested primarily in the space:
\[
\fg_{\ol{1}}^{ss}=\{x\in\fg_{\ol{1}}|\ \operatorname{ad}[x,x]\text{ is semisimple}\}.
\]

\subsubsection{Remark}\label{rmk KM case} When $\fg$ is a Kac-Moody Lie superalgebra, the usual definition of rank naturally extends to all elements in $\fg_{\ol{1}}^{ss}$, and for these superalgebras the description of $\DS_xL(\lambda)$ remains independent of which element is chosen of a given rank, if $L(\lambda)$ is finite-dimensional.  However their actions on the larger category of finite-dimensional $\fg$-modules differs significantly.

\subsection{The case of $\fg=\fq_n$ and arc diagrams} On the other hand, for $\fg=\fq_n$ the behavior of the functor $\DS_x$ on simple modules does not reduce to the case when $x^2=0$; in fact an element $x\in\fg_{\ol{1}}^{ss}$ has rank valued in $\frac{1}{2}\mathbb{N}$, whereas if $x^2=0$ we must have $\rank(x)\in\mathbb{N}$.  Further, there are infinitely many $G_0$-orbits on $\fg_{\ol{1}}^{ss}$ in contrast to the self-commuting cone which has only finitely many orbits. We as yet do not understand $\DS_x$ on simple modules for an arbitrary $x\in\fg_{\ol{1}}^{ss}$; however we know what the possible simple constituents are, and in certain cases we know exactly how the functor behaves.

In order to study the simple constituents of $\DS_xL(\lambda)$, we introduce arc diagrams that have similarities to those used for $\mathfrak{gl},\mathfrak{osp}$, and $\mathfrak{p}$.  In these arc diagrams there are \emph{full arcs}, illustrated with solid lines, and \emph{half-arcs}, illustrated with dotted lines.  The half arcs are exactly those emanating from 0.  Each dominant weight $\lambda$ has an associated arc diagram; for instance if $\lambda=7\epsilon_1+4\epsilon_2+\epsilon_3-\epsilon_6-7\epsilon_7$ then $\Arc(\lambda)$ is given by:
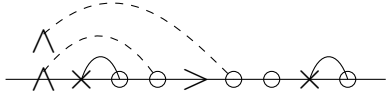
\[
\begin{tikzpicture}
	\draw (0,0) -- (5,0);
	\draw (0.5,0) node[label=center:{\large $\wedge$}] {};
	\draw (0.5,0.5) node[label=center:{\large $\wedge$}] {};
	\draw (1,0) node[label=center:{\large $\times$}] {};
	\draw (1.5,0)  circle(3pt);
	\draw (1,0) .. controls (1.125,0.4) and (1.375,0.4) .. (1.5,0);
	\draw (2,0)  circle(3pt);
	\draw[dashed] (0.5,0.1) .. controls (1, 0.75) and (1.5,0.75) .. (2,0);
	\draw (2.5,0) node[label=center:{\large $>$}] {};
	\draw (3,0)  circle(3pt);
	\draw[dashed] (0.5,0.6) .. controls (1.5, 1.5) and (2.25,0.75) .. (3,0);
	\draw (3.5,0) circle(3pt);
	\draw (4,0) node[label=center:{\large $\times$}] {};
	\draw (4.5,0) circle (3pt);
	\draw (4,0) .. controls (4.125,0.4) and (4.375,0.4) .. (4.5,0);
\end{tikzpicture}
\]
In the above picture, the symbols $\wedge$ lie at 0, and all other symbols lie at positions of $\mathbb{N}_{>0}$.

Write $P^+(\fg)$ for the dominant weights of a Lie superalgebra $\fg$.  The main theorem is as follows:

\subsection{}
\begin{thm}{thmDSqs_intro}
	Take  $\lambda\in P^+(\fq_n)$ and $\nu\in P^+(\fq_{n-2s})$ where $s\in\frac{1}{2}\mathbb{N}$, and let $x\in\fg_{\ol{1}}^{ss}$ be of rank $s$.
	\begin{enumerate}
		\item If $[\DS_x (L(\lambda)):L(\nu)]\neq0$, then $\Arc(\nu)$ can be obtained from $\Arc(\lambda)$ by successively removing $s$ maximal arcs. 
		\item The following are equivalent
\begin{itemize}
\item
 $\operatorname{smult}(\DS_xL(\lambda),L(\nu))\not=0$.
\item
$[\DS_x(L(\lambda)):L(\nu)]=1$;

\item $\zero(\lambda)-\zero(\nu)=2s$; 

\end{itemize}
		\item The indecomposable summands of $\DS_x(L(\lambda))$ are isotypical.
	\end{enumerate}
\end{thm}
Here $\operatorname{smult}(\DS_x(L(\lambda)),L(\nu))=0$ means that $L(\nu)$ appears in $\DS_xL(\lambda)$ the same number of times as $\Pi L(\nu)$; and if $L(\nu)\cong\Pi L(\nu)$ then it further means that $L(\nu)$ appears an even number of times in $\DS_x(L(\lambda))$.  Part (ii) is proven in \cite{GSS}, as well as in section 5.

Note that
for a dominant $\fq_n$-weight $\lambda$ there exists at most one dominant $\fq_{n-2s}$-weight
$\nu$ such that the sets of non-zero coordinates
of $\nu$ and of $\lambda$ coincide; such $\nu$ exists if and only if
$\lambda$ has at least $2s$ of zero coordinates.

\subsection{Case of $\operatorname{rank}x\leq 1$ and $x=C_r$} For $x$ with $\rank x\leq 1$, we are able to precisely describe the multiplicities of the composition factors in a similar manner to the presentations of~\cite{HW},\cite{GH},  see 
Theorems~\ref{main thm rk 1}, \ref{thm rk 1 half}.  
Further, for each $r\leq n$, we may consider the particular element
\[
C_r:=\begin{bmatrix}
	0 & B_r\\ B_r & 0
\end{bmatrix}, \ \text{ where } \ B_r=\begin{bmatrix} I_r & 0\\0 & 0\end{bmatrix}.
\]
Then $\DS_{C_r}$ is completely described on simple modules, and has the special property that for a simple module $L$, $\DS_{C_r}(L)$ is either again simple or is zero. 

\subsection{Properties (a)-(e)}
Let $L$ be a simple finite-dimensional $\fq$-module.  We discuss to what extent the properties from Section \ref{KMp} extend to $\fg=\fq_n$.

\begin{itemize}
\item [$(a)$]
$\DS_x(L)$ is either simple or zero if $\rank x=\frac{1}{2}$; 
the multiplicities of irreducible constituents of $\DS_x(L)$ 
are at most $2$ if $x$ is of rank $1$;

\item[$(b)$] the following are equivalent
\begin{itemize}
\item $L$ is typical;

\item $\DS_x(L)=0$ for all non-zero $x$;

\item 
$\DS_x(L)=0$ for all $x$ of rank $\leq 1$.
\end{itemize}

\item[$(c)$] if the non-graded multiplicity of a simple $\fg_x$-module $L'$ in $\DS_x(L)$   is not $1$, then this non-graded multiplicity is even and $L'$ appears in
$\DS_x(L)$ the same number of times as $\Pi(L)$, see~\Thm{thmDSqs_intro} (ii).

\end{itemize}

\subsubsection{Remark on  $(b)$}\label{remarkb}  If $\DS_x (L(\lambda))=0$ for all $x$ of rank one, it does not imply that $L(\lambda)$ is typical; rather it implies that $\Arc(\lambda)$ has at most one maximal arc, and if it exists it is a half-arc and $\zero(\lambda)=1$.  An example is given by the adjoint representation of $\mathfrak{psq}_{3}$.  The condition that $\DS_x (L(\lambda))=0$ for $x$ of rank one-half is equivalent to asking that $\Arc(\lambda)$ has no half-arcs, in other words $\zero(\lambda)=0$.

\subsubsection{Property $(d)$}
Part $(iii)$ of Theorem \ref{thmDSqs_intro} states that
\begin{itemize}
\item[$(d')$]  $\Ext^1(L',L'')=0$ 
if $L''\not\cong L',\Pi L'$ are simple subquotients 
of $\DS_x(L)$
\end{itemize}
 holds in the $\fq_n$-case. If $\fg$ is a finite-dimensional Kac-Moody superalgebra, this property along with $(c)$ implies $(d)$ (i.e. the semisimplicity of $\DS_x(L)$).  
 
In the $\fq_n$-case we do not yet know if $(d)$ holds in general; however $(d')$ does imply the semisimplicity of $\DS_x(L(\lambda))$ for 
$\lambda$ with no zero coordinates (in particular, for half-integral weights).    Further, in Corollary \ref{DS semisimple} we prove that $DS_x(L)$ is semisimple for any simple finite-dimensional module $L$, when $x$ is of rank 1 with $x^2=0$.  

The following example shows that semisimplicity does not hold for $\fsq_n$: for
a dominant weight $\lambda$ of the form 
$\lambda=\sum_{i=1}^p \lambda_i\vareps_i+\sum_{i=p}^{n-2} \lambda_i\vareps_{i+2}$
with $\sum_{i=1}^{n-2}\lambda_i^{-1}=0$ we have
$L_{\fq_n}(\lambda)=L_{\fsq_n}(\lambda)$. The $\fq_{n-2}$-module
$L_{\fq_{n-2}}(\sum_{i=1}^{n-2} \lambda_i\vareps_i)$ is a submodule of
$\DS_x(L_{\fq_n}(\lambda))$ for $x$ of rank 1 with $x^2=0$, and
this module is a non-splitting extension of two simple
$\fsq_{n-2}$-modules with the same highest weights.

\subsubsection{Property $(e)$}\label{ee}
Such a property fails completely if consider arbitrary $x\in\fg_{\ol{1}}^{ss}$; it is possible it will hold for those $x$ with $x^2=0$.  However we have not computed $\DS_x(L)$ for $x$ of rank greater than $1$.

\subsection{Grading}  Given an element $x\in(\fq_n)_{\ol{1}}$ with $x^2=0$, we may find a semisimple element $h\in(\fq_n)_{\ol{0}}$ satisfying $[h,x]=x$.  In this way we obtain an action of $h$ on $DS_x$ which commutes with $\fg_x$, and thus induces a grading on the $\fg_x$-module $DS_xV$ for a finite-dimensional $\fq_n$-module $V$ (details explained in Sections \ref{section grading qn} and \ref{section grading DS}).  Thus $h$ will act by a scalar on each composition factor of $DS_xV$.    Further it allows us to view $DS_x$ as a functor from $\fg$-modules to $\fg_x\times\mathbb{C}\langle h\rangle$-modules. 

If $x$ is of rank 1 and $L$ is a finite-dimensional simple module, we have computed explicitly the eigenvalues of $h$ on the composition factors of $DS_xL$, and we obtain the following:

\subsection{}\begin{thm}{intro thm mult free}
	If $x$ is of rank $1$ and $L$ is a simple finite-dimensional $\fq_n$-module, then $DS_xL$ is a semisimple, multiplicity-free $\fg_x\times\mathbb{C}\langle h\rangle$-module.
\end{thm}

We note that for a finite-dimensional Kac-Moody Lie superalgebra $\fg$ we also always have such an element $h$ which will acts on $DS_x$ and commutes with $\fg_x$, giving a grading.  We expect that in these cases we also have that if $L$ is a simple finite-dimensional $\fg$-module, then $DS_xL$ will be multiplicity-free as a module over $\fg_x\times\mathbb{C}\langle h\rangle$.  For $\fgl(m|n)$ the grading on $\DS_1(L)$ was computed in~\cite{HW}.

\subsection{Projectivity} In Corollary \ref{cora} we show that if $L$ is a simple finite-dimensional module, then it is projective if and only if $DS_xL=0$ for all $x$ with $\rank x\leq 1$ (in fact a slightly stronger statement is true).  It would be interesting to understand whether this generalizes to all finite-dimensional modules.

\subsection{Methods}
The approach to computing composition factors and multiplicities of $\DS_xL(\lambda)$ is the same as in the $\osp$-case
(see~\cite{GH}): using suitable translation functors the problem 
is reduced to the case of $\fq_{2s}$, where $x$ is of rank $s$.  For $s$ of rank $\leq 1$, we may do the calculations on $\fq_1$ or $\fq_2$ where they are easily performed.  The case of $x=C_r$ is done by using the $\fg_{\ol{0}}$-invariance of $x$ when $n=r$.

The formula  ($d'$) follows from~\Thm{thmDSqs} (i) and the 
 following fact obtained in~\cite{Gdexnew}:
if $\Ext^1(L(\nu),L(\nu'))\not=0$ for  some distinct dominant weights
$\nu,\nu'$, then the diagram of one for these weights
is obtained from the diagram of  the other one by moving a symbol
$\times$ along one of the arches (where $\wedge$ is considered as a half
of $\times$).
 The same reasoning works
in the Kac-Moody case. Note that the  previous proofs of ($d'$) in the Kac-Moody case
were based on a stronger result:
the existence bipartition of $\Ext^1$-graph compatible with the action
of $\DS$-functor. The latter result does not hold in the $\fq_n$-case.

Our proof of complete reducibility for $DS_x(L)$ when $x$ is of rank $1$ with $x^2=0$ relies on the computation of the grading, and ultimately Theorem \ref{intro thm mult free}.  The grading is computed by reducing to the case of $\fq_2$ as above, using the same techniques.

\subsection{Outlook for higher rank elements}  One of the central unanswered questions that we leave to future work is the computation of composition factors and multiplicities of $\DS_x (L(\lambda))$ for $x$ of rank bigger than $1$.  Unlike in the Kac-Moody cases, where one only has to look at those $x$ with $[x,x]=0$ (see Remark \ref{rmk KM case}), for $\fq_n$ (as previously mentioned) we may obtain very different results for elements in $\fg_{\ol{1}}^{ss}$ lying in distinct $G_0=GL(n)$-orbits (see Section \ref{DSqident} for a description of these orbits).

The $G_0$-orbits are parametrized by two pieces of data: first is the rank, a half-integer $r$.  Once we fix a rank $r$, 
the orbits become parametrized by $2r$ unordered complex numbers, where an even number of them are 0 (they are given by the eigenvalues of the matrix $B$ in \ref{intro example}).  This space is an open subvariety of $\mathbb{C}^{2r}/S_{2r}$.  It follows that the behavior of $\DS_{x}$ for $x$ of rank $r$ will be determined by values of symmetric functions on the $2r$ complex numbers; or more precisely, due to the canonical equality $\DS_{cx}=\DS_{x}$ for all $c\in\mathbb{C}^{\times}$, the behavior of $\DS_{x}$ should be determined by the \emph{zero sets} of homogeneous symmetric functions on the $2r$ complex numbers; for an example see the computation of $\DS_x(\fpsq_n)$ below. 

\subsection{Kac-Wakimoto Conjecture and depth} It is of interest to understand to what extent a version of the (generalized) Kac-Wakimoto conjecture, proven in \cite{Skw} for Kac-Moody Lie superalgebras, holds for $\fq_n$.  Recall that the results of \cite{Skw} show that in the Kac-Moody case  we have $\DS_x(L)\neq 0$ if 
$L$ is a simple finite-dimensional module of atypicality\footnote{For $\fq_n$ there are several notions of atypicality, see for instance~\cite{Br} and~\cite{Gcore}.} $r$ and 
$\rank x\leq r$. 

For $\fg=\fq_n$ such a result cannot hold (see~\Rem{remarkb}). Nevertheless one may ask, given $\lambda$ of atypicality $r$, does there exist $x$ of rank $r$ satisfying $\DS_x(L(\lambda))\neq 0$.  

Another possible formulation is the following: given $\lambda$ of atypicality $r$, does there exist a chain of $\DS$ functors, $DS_{x_1},\dots,DS_{x_k}$, such that $\sum\rank(x_i)=r$ and 
\[
(DS_{x_k}\circ\dots\circ DS_{x_1})L(\lambda)\neq0.
\]
We give a positive answer to this question via the notion depth, introduced in Section \ref{depth} (note the definition differs slightly from the one in~\cite{Gcore}). We show that, as in the Kac-Moody case, 
$\depth(N)$ does not exceed the atypicality of $N$ and these numbers are equal if 
$N$ is a simple finite-dimensional module.

\subsection{$\DS_x$ for queer-type superalgebras}\label{intro example} As an illustration of these ideas, we have computed the value of the functor $\DS_{x}$ on all queer-type algebras, i.e. $\fq_n$, $\fsq_n$, $\fpq_n$, and $\fpsq_n$.  In each case we obtain a new superalgebra, and we compute its structure. 
If the rank of $x$ is less than $n/2$, then $\DS_x(\fq_n)=\fq_{n-2r}$
with the similar formulae for  $\fsq_n$, $\fpq_n$, and $\fpsq_n$.
For  $x$ of rank $n/2$ we have $\DS_x(\fq_n)=0$, $\DS_x(\fsq_n)=\mathbb{C}$ and
$\DS_x(\fpq_n)=\Pi\mathbb{C}$. For $\fpsq_n$  we have the following interesting behavior: taking $x=\begin{pmatrix} 0 & B\\B &0\end{pmatrix}$ of rank $n/2$ and letting  $c_1,\dots,c_n$ be the eigenvalues of $B$ we obtain
$$\DS_x(\fpsq_n)\cong\ \left\{
\begin{array}{lcl}
\mathbb{C}^{1|1} & \text{ if } & \sum\limits\limits_i c_{1}\cdots \widehat{c_i}\cdots c_n=\sum\limits_{i}c_{1}^3\cdots \widehat{c_i^3}\cdots c_n^3=0\\
\fq_1  & \text{ if } & \sum\limits\limits_i c_{1}\cdots \widehat{c_i}\cdots c_n=0\neq \sum\limits_{i}c_{1}^3\cdots \widehat{c_i^3}\cdots c_n^3\\
0 & & \sum\limits\limits_i c_{1}\cdots \widehat{c_i}\cdots c_n\neq 0
\end{array}
\right.$$
where by $\mathbb{C}^{1|1}$ we mean the $(1|1)$-dimensional abelian Lie superalgebra (and $\widehat{(-)}$ stands for the exclusion).

\subsection{Acknowledgments}  The authors are grateful to D.~Grantcharov, 
N.~Grantcharov, \\ T.~Heidersdorf, ~V.~Hinich, and V.~Serganova for numerous 
helpful discussions.   The first author was supported by ISF Grant 1957/21.  The second author was supported by ISF Grant 711/18 and NSF-BSF Grant 2019694.

 \subsection{Index of  frequently used notation} \label{sec:app-index}
Throughout the paper the ground field is $\mathbb{C}$; 
$\mathbb{N}$ stands 
 for the set of non-negative integers.  We will frequently used the following notation.

\begin{center}
\begin{tabular}{lcl}
$\fg_{\ol{1}}^{ss}$ & & \ref{gss} \\
$T_{A,B},\ \fh,\ \ft, h_i,\ H_i, \zero(\lambda),\ \nonzero(\lambda)$ & & \ref{Qnotation} \\
$\Cl(\lambda),\ B_{\lambda},\ C_{\lambda},\ L(\lambda) $ & & \ref{Clifford} \\
$\atyp,\ \ \core $, core-free & & \ref{core}\\
$\fq(J),\ \Delta(J)$ & & \ref{notatqJ}\\
$x_{ss}, x_{nil}$ & & \ref{DSqident} \\
$\howl$ & & \ref{howl}\\
$\depth$ & & \ref{depth}\\
$\operatorname{smult}$ & & \ref{subsec smult}\\
\end{tabular}
\end{center}

\section{Preliminaries}\label{prelim}
Throughout the paper the ground field is $\mathbb{C}$; 
$\mathbb{N}$ stands 
 for the set of non-negative integers. We denote by $\Pi$ the parity change functor. 
In Sections~\ref{prelim}--\ref{verblud} 
we  ``identify'' the modules
$N$ and $\Pi(N)$ (where $\Pi$ stands for the parity change functor).
For a finite length module $N$ and a simple module $L$ 
we denote by $[N:L]$ the ``non-graded multiplicity'' (the number of simple subquotients
in a Jordan-H\"older series of $N$ which are isomorphic to $L$ or to $\Pi L$.
We say that a module $N$ is  {\em isotypical} if $N$ is indecomposable and
for each  simple subquotient $L,L'$ of $N$ one has $L'\cong L$ (up to the parity change).

We denote by $\Fin(\fg)$ the  full subcategory of finite-dimensional modules
which are semisimple over $\fg_{\ol{0}}$.

\subsection{$Q$-type superalgebras}\label{Qnotation}
In this paper $\fg$ is the queer ($Q$-type) Lie superalgebra
$\fq_n$. Recall that
 $\fq_n$ is a subalgebra of $\fgl(n|n)$ consisting of the matrices with the block form
$$T_{A,B}:=\begin{pmatrix}
A & B\\
B & A
\end{pmatrix}$$

One has $\fg_{\ol{0}}=\fgl_n$. The group $GL_n$ acts on $\fg$ by the inner
automorphisms;
all triangular decompositions of $\fq_n$ are $GL_n$-conjugated.
We denote by $\ft$ the  Cartan subalgebra of $\fgl_n$
spanned by the elements $h_i=T_{E_{ii},0}$ for $i=1,\ldots,n$.
 Let  $\{\vareps_i\}_{i=1}^n\subset\ft^*$ 
 be the   basis dual to $\{h_i\}_{i=1}^n$.  The algebra $\fh:=\fq_n^{\ft}$ is
a Cartan subalgebra of $\fq_n$; one has $\fh_{\ol{0}}=\ft$.
The elements $H_i:=T_{0,E_{ii}}$ form a basis of $\fh_{\ol{1}}$; one has
$[H_i,H_j]=2\delta_{ij}h_i$.

We fix a usual triangular decomposition:
$\fg=\fn^-\oplus\fh\oplus \fn$,
where $\Delta^+=\{\vareps_i-\vareps_j\}_{1\leq i<j\leq n}$.

We write  $\lambda=\sum_{i=1}^n \lambda_i\vareps_i\in\ft^*$
as  $\lambda=(\lambda_1,\ldots,\lambda_n)$ and  set
$$\zero(\lambda):=\#\{i|\ \lambda_i=0\},\ \ \ \ \nonzero(\lambda):=\#\{i|\ \lambda_i\not=0\}.$$

\subsection{Modules $C_{\lambda}$ and $L(\lambda)$} \label{Clifford} 
The algebra $\cU(\fh)$ can be naturally viewed  as a Clifford algebra over the polynomial algebra $\cS(\ft)$
with the symmetric bilinear form 
$$B:\fh_{\ol{1}}\otimes \fh_{\ol{1}}\to\cS(\ft)\ \ \text{ with } \ \ \ 
B(H,H')=[H,H'].$$

For each $\lambda\in\ft^*$ 
the evaluation of $B$ gives the symmetric form 
$B_{\lambda}: (H,H')\mapsto\lambda([H,H'])$ which defines the Clifford algebra
$$
\Cl(\lambda):=\Cl(\fh_{\ol{1}},B_{\lambda})=\cU(\fh)/\cU(\fh) I(\lambda)$$
 where
$I(\lambda)$ stands for the kernel
of the algebra homomorphism $\cS(\ft)\to\mathbb{C}$ induced by $\lambda$. 
Since $(H_i,H_j)=2\delta_{ij} h_i$ one has
$$\rank B_{\lambda}=\nonzero (\lambda).$$
If $\zero(\lambda)=0$, then each finite-dimensional
 $\Cl(\lambda)$-module is projective.

The algebra $\Cl(\lambda)$ admits a unique simple module  $C_{\lambda}$ (up to isomorphism
and parity change). One has
$\dim C_{\lambda}=2^{[\frac{\rank B_{\lambda}+1}{2}]}$
(and $\sdim C_{\lambda}=0$ for $\lambda\not=0$).  

We denote by $L_{\fg}(\lambda)$ a simple $\fg$-module of the highest weight $\lambda\in\ft^*$;
this module is a unique simple quotient of $\Ind^{\fg}_{\fb} C_{\lambda}$, where 
$C_{\lambda}$  is viewed as a $\fh+\fn^+$-module with the zero action of $\fn^+$.
One has $L_{\fg}(\lambda)_{\lambda}=C_{\lambda}$. We will often omit the index $\fg$
(resp., $\fg_x$) in the notation $L_{\fg}(\lambda)$ (resp., $L_{\fg_x}(\nu)$.

\subsection{Core}\label{core}
The weight $\lambda$ is typical if  $\lambda_i+\lambda_j\not=0$
for all $i,j$; in particular we require that $\lambda_i\neq0$ for all $i$.  We define $\atyp\lambda$ to be
the number of nonzero disjoint pairs $\lambda_i+\lambda_j=0$
for $i\not=j$, plus $\zero(\lambda)/2$.  For example, if $\lambda=(2,1,0,0,0,-2,-3)$, then $\atyp\lambda=5/2$.

Let $\core\lambda$ be
the set obtained from $\{\lambda_i\}_{i=1}^n$ by erasing all zeroes from $\lambda$ and all pairs
$\lambda_i,\lambda_j$ with $i\not=j, \lambda_i+\lambda_j=0$. 
For instance $\core(3,2,0,0,0,-2,-2)=\{3,-2\}$ and this weight has atypicality 
$5/2$. We say that $\lambda$ is {\em core-free} if $\core(\lambda)=\emptyset$.

We denote by $\chi_{\lambda}$ the central character of $L(\lambda)$.
By~\cite{Serq}, 
$\chi_{\lambda}=\chi_{\nu}$
if and only if $\core\lambda=\core\nu$. In particular, the core-free weights are 
the weights with the central character equal to $\chi_0$.
We set
$\atyp\chi_{\lambda}:=\atyp\lambda$. 

\subsection{Subalgebras $\fq(J)$}\label{notatqJ}
Set $I_n:=\{1,\ldots,n\}$.
For each $J\subset  I_n$ we denote by $\fh(J)$ (resp., $\ft(J)$)  the  span
of $H_i,h_i$ (resp., of $h_i$) with $i\in J$ and set
$$\Delta(J):=\{\vareps_i-\vareps_j|\ i,j\in J, \  i\not=j\}.$$
We denote by  $\fq(J)$ the subalgebra spanned by 
$\fh(J)$ and $\sum_{\alpha\in\Delta(J)} \fg_{\alpha}$.
Clearly, $\fq(J)\cong \fq_{k}$, where $k$ is the cardinality of
$J$;  we set 
$\fh_x:=\fh(I_n\setminus J)$
and $\ft_x:=(\fh_x)_{\ol{0}}$ ($\ft_x$ is spanned by 
$h_i$ with $i\in I_n\setminus J$).  The triangular decomposition of
$\fg$ induces  a triangular
decomposition  of $\fq(J)$ with $\Delta(\fq(J))^+=\{\vareps_i-\vareps_j| \ i<j, \ i,j\in J\}$. We will use analogous notation $\fsq(J),\fpq(J),\fpsq(J)$.

\section{$\DS$-functor in $\fq_n$-case}\label{DSinqn}	
The $\DS$-functor was introduced in \cite{DS}.  We recall definitions
and some results of~\cite{DS},\cite{GHSS} in the Appendix.
A study of the $\DS$-functor for $\fq_n$ when $x^2=0$ was initiated in~\cite{S}, see also~\cite{Gcore}, Section 5 for the proofs. In this article we will consider arbitrary $x$ with $x^2$ semisimple; this means $x$ is of the form $x=T_{0,B}$, where $B^2$ is semisimple.  If we write $B=B_{ss}+B_{nil}$ for the Jordan-Chevalley decomposition of $B$, we have $B_{nil}^2=0$ and $B_{ss}B_{nil}=0$.  Define $\rank x=\rank B_{nil}+\rank B_{ss}/2$; in~\Prop{DSg} we show that  $\DS_x (\fq_n)=\fq_{n-2\rank x}$
as well as compute $\DS_x(\fg)$ for other $Q$-type superalgebras.

\subsection{Representatives of $G_0$ on $\fg_{\ol{1}}^{ss}$}\label{DSqident}  
Retain notation of~\ref{notatqJ}.
Recall that the action of $G_0=GL(n)$ on $\fg_{\ol{1}}^{ss}$ is given by the adjoint action on $n\times n$-matrices.  Therefore, finding the $G_0$-orbits on $\fg_{\ol{1}}^{ss}$ is equivalent to finding the $GL(n)$ orbits on square-semisimple $n\times n$ matrices. We give below (non-unique) representatives of these orbits which are easy to work with for our purposes. 
Choose a subset $J:=\{i_p\}_{p=1}^{r}\subset I_n$. Let $s\in\mathbb{N}$ with $s\leq r/2$, and for $p=1,\ldots, s$ fix
a non-zero odd element $x_p\in \fg_{\vareps_{i_{2p-1}}-\vareps_{i_{2p}}}$.  Further, let $c_{2s+1},\dots,c_{r}\in\mathbb{C}^{\times}$, and set
\begin{equation}\label{xform}
x=x_{nil}+x_{ss}:=\sum_{p=1}^s x_p+\sum\limits_{j=2s+1}^{r}c_jH_{i_j}
\end{equation}
Then $x\in \fg_{\ol{1}}^{ss}$ and $x$ has rank $r/2$.

We will always choose $x\in\fg_{\ol{1}}^{ss}$ of the above form and we will use
the above identification of $\fg_x$ with $\fq(I_n\setminus J)$.  Further we will say in this case that $x$ corresponds to the set $J$.  In Section \ref{findimq} we will show in what sense $\DS_x$ is independent of the subset $J$.  

The following lemma is an immediate consequence of Lemma \ref{lemdsxy}.
\subsubsection{}\begin{lem}{ds dim bound lemma}
	Let $x=x_{ss}+x_{nil}$, and write $\ol{x}_{ss}$ (resp., $\ol{x}_{nil}$) for the image of $x_{ss}$ (resp. $x_{nil}$) in $\fg_{x_{nil}}$ (resp., $\fg_{x_{ss}}$).  Let $N$ a finite-dimensional $\fg$-module with semisimple action of $\fg_{\ol{0}}$.  Then
	\[
	\dim \DS_{x}(N)\leq \min\bigg(\dim \DS_{\ol{x}_{ss}}\circ \DS_{x_{nil}}(N),\dim \DS_{\ol{x}_{nil}}\circ \DS_{x_{ss}}(N)\bigg).
	\]
	Further, if $x_{nil}$ is of rank $s$ and we write $x_{nil}=x_1+\dots+x_s$, where each $x_i$ is a rank $1$ nilpotent operator with $[x_i,x_j]=[x_i,x_{ss}]=0$, then for each $0\leq r\leq s$ we have
	\[
	\dim \DS_{x}(N)\leq \dim \DS_{\ol{x}_1}\cdots \DS_{\ol{x}_r}\circ \DS_{\ol{x}_{ss}}\circ\DS_{\ol{x}_{r+1}}\circ\cdots\circ \DS_{x_s}(N),
	\]
	where $\ol{x_i}$ denotes the projection to the appropriate subquotient.	
\end{lem}

\subsection{}
\begin{prop}{DSg}
Take $x$ of rank $r/2$ corresponding to $J\subseteq I_n$ as in~\ref{DSqident}. 
\begin{enumerate}
\item
$\DS_x(\fq_n)$ can be identified with
$\fq(I_n\setminus J)\cong\fq_{n-r}$;

\item $\DS_x$ maps the standard $\fq_n$-module, $L(\vareps_1)$, to the
 standard $\fq_{n-r}$-module;

\item
if $r<n$, then $\DS_x(\fsq_n)$, $\DS_x(\fpq_n)$, $\DS_x(\fpsq_n)$  can be identified with
$\fsq(I_n\setminus J)$, $\fpq(I_n\setminus J)$, $\fpsq(I_n\setminus J)$ respectively;

\item if $r=n$, then $\DS_x(\fsq_n)\cong\mathbb{C}$, $\DS_x(\fpq_n)\cong\Pi\mathbb{C}$;

\item if $r=n>1$, then $\DS_x(\fpsq_n)$ is a commutative $(1|1)$-dimensional Lie superalgebra
if $x=x_{nil}$ or if $\sum c_i^{-1}=\sum c_i^{-3}=0$, $\DS_x(\fpsq_n)\cong \fq_1$
if $\sum c_i^{-1}=0\not=\sum c_i^{-3}$ and
$\DS_x(\fpsq_n)=0$  if $\sum c_i^{-1}\not=0$.
\end{enumerate}
\end{prop}

\begin{proof}
All formulae  can be easily checked in the following cases:
$x=x_{ss}=\sum_{p=1}^r c_i H_i$ or
 $x^2=0$ and $\rank x=1$ (in this case $J=\{i_1,i_2\}$).
Now consider the general case.  It is easy to see that 
$[x,\fq(I_n\setminus J)]=0$ and that $[x,\fq_n]\cap \fq(I_n\setminus J)=0$. Therefore
$\fq(I_n\setminus J)$ with a subalgebra of $\DS_x(\fq_n)$.
Similarly, the standard $\fq(I_n\setminus J)$-module is a submodule
of $\DS_x(L(\vareps_1))$ and, if $r<n$, then $\fsq(I_n\setminus J)$, $\fpq(I_n\setminus J)$, $\fpsq(I_n\setminus J)$  are subalgebras of $\DS_x(\fsq_n)$, $\DS_x(\fpq_n)$, $\DS_x(\fpsq_n)$ respectively.

Using~\Lem{ds dim bound lemma}  we obtain
\begin{equation}\label{DSx1xp}
\dim \DS_x(N)\leq \dim \DS_{\ol{x}_1}\circ \DS_{\ol{x}_2}\circ\ldots\circ \DS_{\ol{x}_p}\bigl(\DS_{x_{ss}}(N)\bigr)
\end{equation}
for any finite-dimensional module $N$. Using the cases $x=x_{ss}$ and $\rank x=1$
we get $\dim \DS_x(\fq_n))\leq \dim \fq(I_n\setminus J)$ and the similar inequalities
for other cases. This establishes (i), (ii), (iii). 

Consider the remaining case $\rank x=n/2$.
By (i), $\DS_x(\fq_n)=0$.
Using Hinich's lemma we get (iv).  

For (v) we  have $J=\{1,\ldots,n\}$; thus $x$ takes the form
\[
x:=\sum_{p=1}^s x_p+\sum_{j=2s+1}^n c_j H_{j},
\]
with  $c_j\in\mathbb{C}^\times$, $x_p\in\fg_{\vareps_{2p-1}-\vareps_{2p}}$.  We set
\[
\ \ 
y:=\sum_{p=1}^s y_p+\sum_{p=2s+1}^nc_i^{-1}H_{i_p},
\]
where $y_p\in \fg_{\vareps_{2p}-\vareps_{2p-1}}$ is an odd element satisfying $[x_p,y_p]=2h_{2p-1}+2h_{2p}$.  In particular, we have $[x,y]=2T_{\id,0}$.

Applying Hinich's Lemma to the short exact sequence
\[
0\to \mathbb{C}T_{\Id,0}\to \fsq_n\to \fpsq_n\to 0
\]
and using $\DS_x(\fsq_n) =\mathbb{C}$ 
we obtain a long exact sequence
\begin{equation}\label{long}
0\to E\to \mathbb{C}\to \mathbb{C}\to \DS_x(\fpsq_n)\to \Pi E\to 0
\end{equation}
where $E:=\mathbb{C}T_{\Id,0}\cap [x,\fsq_n]$. It is easy to see that
 $E=\mathbb{C}T_{\Id,0}$ if $y\in\fsq_n$
and $E=0$ otherwise. This gives $\dim\DS_x(\fpsq_n)=(1|1)$ 
if $x=x_{nil}$ or if $\sum c_i^{-1}=0$, and
$\DS_x(\fpsq_n)=0$  if $\sum c_i^{-1}\not=0$.

Consider the case when  $x=x_{nil}$ or  $\sum c_j^{-1}=0$, so that $\dim \DS_x(\fpsq_n)=(1|1)$.
By~(\ref{long}),
$\DS_x(\fpsq_n)_{\ol{0}}$ is an ideal in $\DS_x(\fpsq_n)$, so 
$\DS_x(\fpsq_n)$ is either commutative or isomorphic to $\fq_1$.
It is easy to see that 
$y$ has a non-zero image in $\DS_x(\fpsq_n)$; we denote this image by $\ol{y}$.
Then $\DS_x(\fpsq_n)$ is commutative  if and only if $\ol{y}^2=0$ that is
$y^2\in \mathbb{C}T_{\Id,0}+[\fsq_n,x]$.

If $x=x_{nil}$, then $y^2=0$. 
In the remaining case $\sum c_j^{-1}=0$ and  $y^2=\sum_{j=2s+1}^n c_j^{-2}h_j$. The elements in $[\fsq_n,x]\cap \ft$ 
are of the form $\sum d_ih_i$, where 
$d_{2p-1}=d_{2p}$ for $p=1,\ldots,s$ and $\sum_{j=2s+1}^n c_j^{-1}d_j=0$.
Therefore $y^2\in \mathbb{C}T_{\Id,0}+[\fsq_n,x]$ is equivalent to
the existence of $d\in\mathbb{C}$ and $z_j\in\mathbb{C}$ for $j=2s+1,\ldots,n$ such that
$c_j z_j=c_j^{-2}+d$ for all $j$ and $\sum z_j=0$.
Since $\sum c_j^{-1}=0$, the last formula is equivalent to $\sum c_j^{-3}=0$
This completes the proof.
\end{proof}

\subsection{}
Let $x\in\fg_{\ol{1}}^{ss}$, and write $Z(\fg)$ for the center of the enveloping algebra $\mathcal{U}(\fg)$.  Then we have
\[
Z(\fg)\subseteq \mathcal{U}(\fg)^x\to\mathcal{U}(\fg_x).
\]
Thus we have a natural map $Z(\fg)\to Z(\fg_x)$, inducing a pullback map on central characters
\[
\eta_x^*:\Hom(Z(\fg_x),\mathbb{C})\to\Hom(Z(\fg),\mathbb{C}).
\]
As in $\fgl(m|n)$ and $\osp(m|n)$-cases, the $\DS$-functor respects preserves cores, in the following sense. 

\subsubsection{}\begin{lem}{ss pullback map}  Let $x\in\fg_{\ol{1}}^{ss}$ with $r=\rank x$.  Then the map 
	\[
	\eta_{x}:Z(\fq_{n})\to Z(\fq_{n-2r})
	\]
	is surjective.  Further, we have
$	\eta_x^*(\chi_c)=\chi_{c}$.
\end{lem}
\begin{proof}
	See {\em Theorem} (\cite{S}, Thm. 6.3, \cite{Gcore}, Cor. 5.8.1) for proofs in the case when $x^2=0$; the case of $x\in\fg_{\ol{1}}^{ss}$ is almost identical.  
\end{proof}

\subsubsection{}\begin{cor}{corchi}
If a $\fg$-module $N$ has a central character $\chi_{\lambda}$,
then a $\fg_x$-module $\DS_x(N)$ has a central character $\chi_{\nu}$,
where $\core(\lambda)=\core(\nu)$. 

In particular,
$\atyp\lambda=\atyp\nu-\rank x$ and
 $\DS_x(N)=0$
if $\atyp\lambda<\rank x$.
\end{cor}

\subsection{Independence of $\DS_x$ from $J\subseteq\{1,\dots,n\}$}\label{findimq}
Retain notation of~\ref{DSqident}. Let $I_s:=\{n,n-1\ldots,n-s+1\}$ and let $x_s\in\fg_{\ol{1}}^{ss}$ correspond to $I_s$ as in Section \ref{DSqident}.  Let $J=\{j_1<\cdots<j_s\}\subseteq\{1,\dots,n\}$ be an arbitrary subset of size $s$ and let $\sigma\in S_n\subseteq GL(n)$ denote the permutation with $\sigma(n-s+1)=j_1,\dots,\sigma(n)=j_s$, and set $x:=\sigma(x_s)$.
The action of $\sigma$ gives an isomorphism $\sigma_x:\fg_{x_s}\iso \fg_{x}$ 
and the commutative diagram
$$\xymatrix{& \Fin(\fg)\ar^{\phi}[r] \ar^{\DS_{x_s}}[d] & \Fin(\fg)\ar^{\DS_{x}}[d]& \\
 & \Fin(\fg_{x_s})\ar^{\phi_x}[r]  & \Fin(\fg_{x}) &} $$
where the functors $\phi$ and $\phi_x$ correspond to the shift of module structure
along $\sigma:\fg\to\fg$ and $\sigma:\fg_{x_s}\to\fg_{x}$ respectively.
If $N$ is a finite-dimensional $\fg$-module, the action of $g$  on
$N$ induces an isomorphism $\phi(N)\cong N$ and a bijection between $\DS_{x_s}(N)$ and $\DS_{x}(N)$ which is compatible with the algebra isomorphism $\fg_{x_s}\iso \fg_{x}$.

By~\Prop{DSg}, we have natural identifications $\fg_{x_s}\cong\fq_{n-s}$ and  $\fg_{x}\cong\fq_{n-s}$ as the subalgebras of $\fq_n$ corresponding to the subsets $I_s$ and $J$.  Under these identifications we have a commutative diagram
\[
\xymatrix{\fg_{x_s} \ar[r]^{\sigma} \ar[d]^{\sim} & \fg_x\ar[d]^{\sim} \\ 
	\fq_{n-s} \ar[r]^{\sigma} & \fq_{n-s} }
\]
Further, $\sigma$ takes the Cartan subalgebras and subalgebras of the respective copies of $\fq_{n-s}$ to one another.  It follows that if $\nu\in P^+(\fq_{n-s})$ then $\phi_{x}(L_{\fg_{x_s}}(\nu))\cong L_{\fq_{x}}(\nu)$.  

It follows that the functor $\DS_x: \Fin(\fg)\to\Fin(\fg_x)$ 
is independent of the choice of set $J\subseteq I_n$ that $x$ corresponds to. (In the computations below 
our choice of $x$ will depend on the central character.)

\subsubsection{}\begin{prop}{ss_kills_half_integral}
Let $M$ be a finite-dimensional $\fg$-module such $\zero(\nu)=0$ for any weight
$\nu$ of $M$. If 
 $x\in(\fq_{n})_{\ol{1}}^{ss}$ with $x^2\neq 0$, then $\DS_x (M)=0$.
\end{prop}

\begin{proof}
As in Section \ref{DSqident} we may write $x=x_{ss}+x_{nil}$ 
with $x_{ss}\in\fh$ and $x_{ss}^2=x^2$.  	Set $N:=M^{x^2}$.
There exists $t\in\ft$ such that $[t,x_{nil}]=x_{nil}$.
Since $x_{ss}\in\fh$ we have $[t,x_{ss}]=0$. Since $\zero(\nu)=0$ for any weight
$\nu$ of $M$, 
 by~\ref{Clifford},  $M$ is a projective $\fh$-module, implying that $\DS_{x_{ss}} M=0$ (since $x_{ss}\in\fh$), i.e.  $M^{x_{ss}}=x_{ss} M$.
 Now the assertion follows from~\Lem{lemdsxy}.
 \end{proof}

\subsection{Gradings for elements with $x^2=0$}\label{section grading qn}  
Fix $x\in\fq(J)$ with $x^2=0$ (i.e., $x=x_{nil}$). In this case there exists an element $h\in\ft\cap \fq(J)$ such that  $[h,x]=x$.  Thus, as is explained in Section \ref{section grading DS}, we obtain a grading on $M_x$ as a $\fg_x$-module for any finite-dimensional $\fg$-module $M$. Note that $[h,\fg_x]=0$.

%
%
%
%
%
%
%

\section{Dominant weights and arc diagrams}
In this paper we study the action of $\DS_x$ on finite-dimensional simple modules.
 We denote by $P^+(\fg)$ the set of dominant weights, i.e.
\[
P^+(\fg):=\{\lambda\in\ft^*|\ \dim L(\lambda)<\infty\}.
\]

By~\cite{P}, $\lambda\in P^+(\fg)$ if and only if
 $\lambda_i-\lambda_{i+1}\in\mathbb{N}_{\geq 0}$ and $\lambda_i=\lambda_{i+1}$
implies $\lambda_i=0$. We call weight $\lambda$ {\em integral } (resp., {\em half-integral}) 
if $\lambda_i\in\mathbb{Z}$ (resp., $\lambda_i-\frac{1}{2}\in\mathbb{Z}$) for all $i$. If $\lambda\in P^+(\fg)$ is atypical, then $\lambda$ is either integral or half-integral. By~\ref{corchi}, $\DS_x(N)=0$ for each typical module $N$
and $x\not=0$. Without loss of generality we assume 

\begin{center}{\em $x\not=0$ and
 $\lambda$ is integral or half-integral}.\end{center}

By~\Prop{ss_kills_half_integral},
for half-integral weights it is only interesting to consider those $x$ with $x^2=0$.

\subsection{Weight diagrams}
Weight  diagrams
were first defined in \cite{BS4} for $\mathfrak{gl}(m|n)$. The conventions on how to draw these weight diagrams differ; we follow essentially~\cite{GS}.

For $\lambda=(\lambda_1,\ldots,\lambda_n)$ we 
construct the weight diagram $\diag(\lambda)$ as follows:
\begin{itemize}
\item for $s\not=0$ 
we put $>$  (resp., $<$)  at the position $s$ if
there exists $j$ with $\lambda_j=s$ (resp., $\lambda_j=-s$);
\item
we write $\times$ if
the position $s\neq0$ contains $>$ and $<$;

\item  if $\zero(\lambda)\not=0$, 
we put $\wedge^r$, where $r=\zero(\lambda)$; 

\item if $\lambda$ is integral (resp., half-integral)
 we put the empty
sign $\circ$ at each non-occupied position with the coordinate
 in $\mathbb{N}$ (resp., in $\frac{1}{2}+\mathbb{N}$).
\end{itemize}

If $\lambda$  is integral, 
we draw the diagram from the zero position, and for $\lambda$ half-integral we draw from the $\frac{1}{2}$ position.  For instance,
$$\begin{array}{lll}
(4,1,-1,-3,-4)\ \ & &\ \circ\times\circ<\times\circ\circ\ldots\\
(5,2,0,0,0,0,0,-2-,-3)\ \ & &  \wedge^5\ \circ \times <\circ >\circ\ldots\\
(\frac{5}{2},\frac{1}{2},-\frac{1}{2},-\frac{5}{2}) & & \times\circ\times\circ\ldots
\end{array}
$$
where $\ldots$ stands for the infinite sequence of the empty symbols $\circ$.

We obtain a one-to-one correspondence between integral (resp., half-integral) 
dominant weights
and the diagrams containing $n$ symbols $>,<,\wedge$ (where $\times$ considered as 
the union of $<$ and $>$), where each non-zero position contains 
exactly one of the symbols $>,<,\times$ or $\circ$ and the zero position
contains $\wedge^r$ or $\circ$.   Note that  $\atyp\lambda$ is equal to the number of $\times$s in $\diag(\lambda)$ plus half the number of symbols $\wedge$ at 0.  We roughly think of the symbol $\wedge$ as half of a symbol $\times$.

For a weight diagram $f$ denote by $f(i)$ the symbols at the $i$th position.

\subsubsection{Core diagrams}
The symbols $>,<$ are called {\em core symbols}.

A {\em core diagram} is a weight diagram which does not contain symbols $\times$ and $\wedge$.
The core diagram
of $\lambda$ is obtained from the weight diagram of $\lambda$
by erasing all $\times$ and $\wedge$ symbols.

We say that a $\fg$-central character $\chi$ is dominant if 
there exists a finite-dimensional module with this central character.
By above, the central characters are parametrized by $\core(\lambda)$, so
the dominant central characters of atypicality $k>0$ 
are parametrized by the  core diagrams with 
$n-2k$ non-empty symbols. For instance, the central character of the maximal atypicality corresponds to the
empty core diagram; for $\fq_2$ the diagrams of the weights with such central character are $\wedge^2$, $\circ\times$, $\circ\circ\times$ and so on; for $\fq_3$ these diagrams  are $\wedge^3\circ$, $\wedge\times$, $\wedge\circ\times$ and so on.

For a core diagram $f$ we denote by $\chi(f)$ the corresponding central character.
 
\subsubsection{Core-free diagrams}\label{howl}
Note that $\lambda$ is core-free if and only if
the weight diagram of $\lambda$ does not have core symbols. We assign to each diagram $f$ a core-free diagram $\howl(f)$
which is obtained from $f$ by erasing all core symbols. For instance,
$$\begin{array}{lllll}
f=\wedge^4>\circ\times<\times\circ\ldots & &  \core(f)=\circ>\circ\circ <\circ\ldots & &  
\howl(f)=\wedge^4\circ\times\times\circ\ldots\\
g=\wedge^3\circ\times<\times\circ\ldots & & 
\core(g)=\circ\circ\circ<\circ\ldots & & 
 \howl(g)=\wedge^3\circ\times\times\ldots\end{array}$$

For an atypical dominant weight $\lambda$ we denote by $\howl(\lambda)$ the weight corresponding to the diagram $\howl(\diag(\lambda))$.

\subsection{Partial order }
We will consider the standard partial order
on $\ft^*$: 
$$\lambda>\nu\ \text{ if }\lambda-\nu\in\mathbb{N}\Delta^+.$$

\subsubsection{}
\begin{lem}{lemneworder}
For atypical weights $\eta,\mu\in P^+(\fq_n)$ one has
$$\core(\eta)=\core(\mu)\ \ \&\ \ \eta>\mu\ \ \ \Longrightarrow\ \ \
\howl(\eta)\not\leq\howl(\mu).$$
\end{lem}
\begin{proof}
Let $f,g$ be weight diagrams of atypicality $p>0$ 
with $\core f=\core g$, 
and let
 $a_1\geq a_2\geq \ldots \geq a_p$
(resp., $b_1\geq b_2\geq \ldots \geq b_p$) be
 the  coordinates of the symbols
$\times$ in $f$ (resp., in $g$). 
We write $g\succ f$ if $b_j>a_j$ for some $j$ and
$a_i=b_i$ for each $i<j$. For example
$$ \wedge>\circ\times\times\ \ \succ\ \ \wedge^3>\circ\circ\times.$$
It is easy to see that for atypical weights $\mu,\eta\in P^+(\fq_n)$ with
$\core(\mu)=\core(\eta)$ one has
$$\begin{array}{l}
\eta>\mu\ \ \Longrightarrow\ \ \diag(\eta)\succ\diag(\mu)\\
\diag(\eta)\succ\diag(\mu)\ \ \Longleftrightarrow\ \ \diag(\howl(\eta))\succ\diag(\howl(\mu)).\end{array}$$
If $\eta>\mu$ and $\howl(\eta)\leq\howl(\mu)$, then 
$\diag(\eta)\succ\diag(\mu)$ and $\diag(\howl(\mu))\succ\diag(\howl(\eta))$,
a contradiction.
\end{proof}

\subsubsection{Remark}
Note that 
$$\core(\eta)=\core(\mu)\ \ \&\ \ \howl(\eta)>\howl(\mu)\ \ \ \not\Longrightarrow\ \ \
\eta>\mu.$$
For example, for $\lambda:=(4,1,0,0,-4)$, $\nu:=(3,2,1,-2,-3)$ one has 
$\howl(\lambda)=(3,0,0,-3)$ and $\howl(\nu)=(2,1,-1,2)$. In this case
$\lambda\not>\nu$ and
$\howl(\lambda)>\howl(\nu)$.

\subsection{Arc diagrams}\label{free}
A {\em generalized  arc diagram} is the following data: 
\begin{itemize}
\item[$\bullet$]
a weight diagram $f$, where
the symbols $\wedge$ at the zero position are drawn vertically;
\item[$\bullet$]
a collection of non-intersecting arcs of two types:
\begin{itemize}
\item a \emph{full arc}
$\arc(a;b)$  connects
 the symbol $\times$ at the position $a\not=0$ with the empty symbol
at the position $b>a$; full arcs are depicted by solid arcs;
\item a \emph{half arc}
$\arc(0;b)$  connects
 the symbol $\wedge$ at the zero position with an empty symbol
at the position $b$; half arcs are depicted by dashed arcs.
\end{itemize}

\end{itemize}

An empty position  is called {\em free} if it is not an end of an arc.

For $a\neq0$, we call $\arc(a;b)$ a {\em full arc supported at} $a$, and 
$\arc(0;b)$ a {\em half arc supported at} $0$.  In this sense, when later on we talk about a number of arcs, two half arcs will make one arc.  For example, a quantity of 3/2 arcs consists either of a full arc with a half arc, or three half arcs.

A generalized arc diagram is called an {\em arc diagram} if each symbol $\times$ and $\wedge$ is the left end of exactly one arc and there are no free positions under the arcs.
Each weight diagram $f$ admits a unique arc diagram
which we denote by $\Arc(f)$ (as in $\fgl(m|n)$ and $\osp(m|n)$ case we 
construct the arcs successively starting from the rightmost symbol
$\times$). We write $\Arc(\nu)$
for $\Arc(\diag(\nu))$.
We write $\Arc(\nu)\subset\Arc(\lambda)$ if each arc in $\Arc(\nu)$
appears in $\Arc(\lambda)$.
	
%

\subsubsection{Partial order}
We consider a partial order on the set of arcs by saying that one arc is smaller than
another one if the first one is "below" the second one, that is
\[
arc(a;b)>arc(a';b') \ \text{ if and only if } \ a<a'<b'<b.
\]

Since the arcs do not intersect, one has
\[
arc(a;b)>arc(a';b')\ \ \Longleftrightarrow\ \ a<a'<b,\\
\]
and
any two distinct arcs of the form $arc(0,b_1),arc(0,b_2)$ are comparable: either 
$arc(0,b_1)>arc(0,b_2)$  or
$arc(0,b_1)<arc(0,b_2)$.


%
%
%
%

\subsubsection{Remarks}

Notice that a maximal arc can be ``removed'': if we erase the symbol $\times$ (or $\wedge$) on the left end of the arc and the arc itself, we obtain another arc diagram (this does not hold if the arc is not maximal);
if $\Arc(\nu)$ is obtained from $\Arc(\lambda)$ in this way, then
$\Arc(\nu)\subset\Arc(\lambda)$.

Observe there is the natural one-to-one correspondence
between $\Arc(f)$ and $\Arc(\howl(f))$; this correspondence preserves
the partial order on arcs and the type of arc (full or half).

\subsubsection{Examples}  The arc diagram for the weight $\lambda=7\epsilon_1+4\epsilon_2+\epsilon_3-\epsilon_6-7\epsilon_7$ looks as follows:
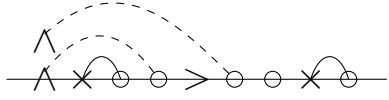
\[
\begin{tikzpicture}
	\draw (0,0) -- (5,0);
	\draw (0.5,0) node[label=center:{\large $\wedge$}] {};
	\draw (0.5,0.5) node[label=center:{\large $\wedge$}] {};
	\draw (1,0) node[label=center:{\large $\times$}] {};
	\draw (1.5,0)  circle(3pt);
	\draw (1,0) .. controls (1.125,0.4) and (1.375,0.4) .. (1.5,0);
	\draw (2,0)  circle(3pt);
	\draw[dashed] (0.5,0.1) .. controls (1, 0.75) and (1.5,0.75) .. (2,0);
	\draw (2.5,0) node[label=center:{\large $>$}] {};
	\draw (3,0)  circle(3pt);
	\draw[dashed] (0.5,0.6) .. controls (1.5, 1.5) and (2.25,0.75) .. (3,0);
	\draw (3.5,0) circle(3pt);
	\draw (4,0) node[label=center:{\large $\times$}] {};
	\draw (4.5,0) circle (3pt);
	\draw (4,0) .. controls (4.125,0.4) and (4.375,0.4) .. (4.5,0);
\end{tikzpicture}
\]
There are two maximal arcs: $\arc(0,5)$ and $\arc(7,8)$.  

The arc diagram for $\lambda=\frac{5}{2}\epsilon_1+\frac{1}{2}\epsilon_2-\frac{1}{2}\epsilon_3-\frac{3}{2}\epsilon_4-\frac{5}{2}\epsilon_5$ is given by:
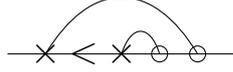
\[
\begin{tikzpicture}
	\draw (0,0) -- (3,0);
	\draw (0.5,0) node[label=center:{\large $\times$}] {};
	\draw (1,0) node[label=center:{\large $<$}] {};
	\draw (1.5,0) node[label=center:{\large $\times$}] {};
	\draw (2,0)  circle(3pt);
	\draw (1.5,0) .. controls (1.625,0.4) and (1.875,0.4) .. (2,0);
	\draw (2.5,0)  circle(3pt);
	\draw (0.5,0) .. controls (1.1, 1) and (1.9,1) .. (2.5,0);
\end{tikzpicture}
\]
In this case there is only one maximal arc: $\arc(\frac{1}{2},\frac{9}{2})$.
\section{Main results}\label{mainresult}
In this section we formulate the main results, illustrate them
by examples and give outlines of the proofs.  In the following, for $x\in\fg_{\ol{1}}^{ss}$ of rank $s$, $\lambda\in P^+(\fq_n)$, and $\nu\in P^+(\fq_{n-2s})$, we write $m_x(\lambda;\nu):=[\DS_x(L_{\fq_n}(\lambda)):L_{\fq_{n-2s}}(\nu)]$.

\subsection{}
\begin{thm}{thmDSqs}
Take  $\lambda\in P^+(\fq_n)$ and $\nu\in P^+(\fq_{n-2s})$ where $s\in\frac{1}{2}\mathbb{N}$, and let $x\in\fg_{\ol{1}}^{ss}$ be of rank $s$.
\begin{enumerate}
\item
If $m_x(\lambda;\nu)\not=0$, then $\core(\lambda)=\core(\nu)$ and
$\Arc(\nu)$ can be obtained from $\Arc(\lambda)$ by successively removing $s$ maximal arcs. 
\item We have $m_x(\lambda,\nu)=1$ if and only if $\zero(\nu)=\zero(\lambda)-2s$.  If $\zero(\nu)\neq\zero(\lambda)-2s$ then
\[
\operatorname{smult}(DS_xL(\lambda),L(\nu))=0.
\]
(See Section \ref{subsec smult} for the definition of $\operatorname{smult}$).
\item Let $x=T_{0,E_{11}+\dots+E_{rr}}$.  Then $\DS_{x}(L(\lambda))=0$ if $\zero(\lambda)<r$, and if $\zero(\lambda)\geq r$ then $\DS_{x}(L(\lambda))=L(\lambda')$, where $\lambda'$ is obtained from $\lambda$ by removing $r$ zeroes; in other words $\Arc(\lambda')$ is obtained from $\Arc(\lambda)$ by removing $r$ half arcs.
\end{enumerate}
\end{thm}

When $x$ is of rank $1/2$ its $G_0$-orbit intersects $\mathbb{C}\langle T_{0,E_{11}}\rangle$, and thus part (ii) of Theorem \ref{thmDSqs} covers this case.  For the case when the rank of $x$ is 1, we have the following more precise statements.

\subsection{}\begin{thm}{main thm rk 1}  Let $x=T_{0,B}\in\fg_{\ol{1}}^{ss}$ be of rank 1 and $\lambda$ be a dominant integral weight.
	\begin{enumerate}
	\item If $\lambda$ is integral and
$\tr B=0$ then $m_x(\lambda;\nu)=2$ if $\Arc(\nu)$ can be obtained from $\Arc(\lambda)$ by removing 
	one maximal full arc,  $m_x(\lambda;\nu)=1$ if $\Arc(\nu)$ can be obtained from $\Arc(\lambda)$ by removing two successive
	maximal half arcs and $m_x(\lambda;\nu)=0$ in other cases.

\item If $\lambda$ is integral and
$\tr B\not=0$ then $m_x(\lambda;\nu)=1$ if $\Arc(\nu)$ can be obtained from $\Arc(\lambda)$ by removing two successive
	maximal half arcs and $m_x(\lambda;\nu)=0$ in other cases.
\end{enumerate}
\end{thm}

Recall that $\DS_x(L(\lambda))=0$ if $\lambda$ is a dominant half-integral weight 
and $x^2\not=0$, see~\Prop{ss_kills_half_integral}.

\subsection{}\begin{thm}{thm rk 1 half}
Let $x\in\fg_{\ol{1}}^{ss}$ be of rank 1 with $x^2=0$. 
 If $\lambda$ is half-integral, then  $m_x(\lambda;\nu)=2$ if $\Arc(\nu)$ can be obtained from $\Arc(\lambda)$ by removing 
	one maximal full arc and $m_x(\lambda;\nu)=0$ in other cases.
\end{thm}

\subsubsection{Remark}
If $\Arc(\nu)$ is obtained from $\Arc(\lambda)$ by removing $r$ maximal arcs and $s$ maximal half arcs, then $\zero(\nu)=\zero(\lambda)-s$.


%
%

\subsubsection{}
\begin{cor}{cora}
For an atypical $\lambda\in P^+(\fq_n)$ the following assertions are equivalent:
\begin{enumerate}
\item
$\DS_x(L(\lambda))=0$ with $x^2=0$, $\rank x=1$; 
\item
$\Arc(\lambda)$ has a unique maximal arc and $\zero(\lambda)=1$;
\item $\DS_x(L(\lambda))=0$ for any $x\not=0$ with $\rank x\in\mathbb{Z}$.
\end{enumerate}
\end{cor}

\subsubsection{}\begin{cor}{}
	For every simple module $L(\lambda)$ with $\atyp\lambda>0$, there exists a nonzero $x\in\fg_{\ol{1}}^{ss}$ of rank $\leq 1$ such that $\DS_xL(\lambda)\neq0$.
\end{cor}
\begin{proof}
	Suppose that $\DS_xL(\lambda)=0$ for all nonzero $x\in\fg_{\ol{1}}^{ss}$ of rank 1; then by Corollary \ref{cora}, $\zero(\lambda)=1$.  Therefore by Theorem \ref{main thm rk 1}, $\DS_{C_1}L(\lambda)\neq0$.
\end{proof}

\subsubsection{Remark} 
It is not hard to see that if  the coordinates of $\nu$ are obtained
from the coordinates of $\lambda$ by removal $2r$ zero coordinates, then 
$m_x(\lambda;\nu)=1$  for $x$ of rank $r$ (see~\cite{Gcore}, Prop. 5.7.2
for slightly more general statement); in this case
$\Arc(\nu)$ is obtained from $\Arc(\lambda)$ by removal of $2r$ half arcs.

\subsection{}
\begin{cor}{corcomp} Let $x\in\fg_{\ol{1}}^{ss}$.  
\begin{enumerate}

\item All indecomposable components of $\DS_x(L_{\fq_n}(\lambda))$
are isotypical.

\item
 If $\zero(\lambda)=0$, then $\DS_x(L(\lambda))$
is completely reducible.
\end{enumerate}\end{cor}

\begin{proof}
In order to prove that each indecomposable component
of $\DS_x(L(\lambda))$ is isotypical, it is enough to verify that 
$\Ext^1(L(\nu),L(\mu))=0$ if $\mu\not=\nu$ and
$L(\mu), L(\nu)$ are subquotients of  $\DS_x(L(\lambda))$.
Take $\mu,\nu$ as above.
By~\Thm{thmDSqs}
all arcs in $\Arc(\mu)$, $\Arc(\nu)$ are also arcs in $\Arc(\lambda)$;
in particular, if $\diag\lambda$ has $\circ$ at some position, then both diagrams
$\diag(\nu),\diag(\mu)$ also have  $\circ$ at this position.
Assume that $\mu\not\geq \nu$ and  that 
$\Ext^1(L(\nu),L(\mu))\not=0$ or $\Ext^1(L(\mu),L(\nu))\not=0$.
By~\cite{Gdex}, Theorem  A, this implies that   $\diag(\nu)$ can be obtained
from  $\diag(\mu)$ by moving one symbol $\times$ in $\diag(\mu)$
along the arc in $\Arc(\mu)$ or changing two symbols $\wedge$ 
to a symbol $\times$ which occupies a position connected to one of the symbols
$\wedge$ in $\Arc(\mu)$. This means that
there are two positions $a<b$ which
are connected by an arc in $\Arc(\mu)$ ($\diag(\mu)$ has $\circ$ at the position $b$) 
and $\diag(\nu)$ has $\times$ at the position $b$.
By above, the positions $a,b$
are connected by an arc in $\Arc(\lambda)$, so
$\diag\lambda$ has $\circ$ at the position $b$ and thus 
$\diag(\nu)$ has $\circ$ at the position $b$, a contradiction.
This establishes  (i). 

If $\zero(\lambda)=0$ and
$L(\nu)$ is a subquotient of  $\DS_x(L(\lambda))$, then
 $\zero(\nu)=0$; by
Theorem 3.1 in~\cite{Nicki} this gives
$\Ext^1(L(\nu), L(\nu))=\Ext^1(L(\nu), \Pi L(\nu))=0$.
Combining with (i) we deduce complete reducibility of $\DS_x(L(\lambda))$.
\end{proof}

\subsection{A result on grading and semisimplicity}
Take $x$ with $x^2=0$ to be of rank $1$.
  
  \subsubsection{Notation}
Let $\lambda\in P^+(\fq_n)$ be integral or half-integral.  For $k\in\frac{1}{2}\mathbb{N}$, we let 
	\[
g(\lambda,k)=\left\{\begin{array}{ll}
	 \ell(\lambda,k)+1 & \text{ if }\lambda \text{ integral}\\
	\ell(\lambda,k)+1/2 & \text{ if }\lambda \text{ half-integral}.\\
\end{array}\right.
\]
Here $\ell(\lambda,k)$ denotes the number of nonzero free positions strictly to the left of $k$ on the arc diagram of $\lambda$ (see~\ref{free} for definition of a free
position).  For example, if $\lambda$ is half-integral and has weight diagram $\times\circ\circ\times$ then $g(\lambda,7/2)=3/2$ and $g(\lambda,1/2)=1/2$, while if $\lambda$ is integral with weight diagram $\wedge^2\times\circ\times$, then $g(\lambda,3)=g(\lambda,1)=1$. 

\subsubsection{}\begin{thm}{thm grading}
Let $x$ be of rank 1 with $x^2=0$, and let $h$ be as in Section \ref{section grading qn}.  Let $\lambda\in P^+(\fq_n)$. 
	\begin{enumerate}
		\item If $\mu\in P^+(\fq_{n-2})$ is obtained from $\lambda$ removal of 2 zeros, then $h$ acts trivially on $L(\mu)$ as a submodule of $\DS_xL(\lambda)$.
		\item If $\mu\in P^+(\fq_{n-2})$ is obtained from $\lambda$ by removal of a symbol $\times$ at position $k>0$ in its arc diagram, then $h$ acts with opposite eigenvalues $\pm g(\lambda,k)$ on the two copies of $L(\mu)$ in $\DS_xL(\lambda)$.
	\end{enumerate}  
\end{thm}

\subsubsection{Remark} Theorems \ref{main thm rk 1} and \ref{thm grading} imply that $DS_xL(\lambda)$ is multiplicity-free as a module over $\fg_x\times\mathbb{C}\langle h\rangle$.  Thus we obtain as a simple corollary:

\subsection{}\begin{cor}{DS semisimple}
	For $x$ of rank $1$ with $x^2=0$, $\DS_xL(\lambda)$ is semisimple.
\end{cor}
\begin{proof}
	This follows immediately from Corollary \ref{corcomp} and Theorem \ref{thm grading}.
\end{proof}

\subsubsection{Example}  If $\lambda=7\epsilon_1+4\epsilon_2+2\epsilon_3-2\epsilon_5-4\epsilon_6-7\epsilon_7$, then the corresponding arc diagram looks as follows:
\[
\begin{tikzpicture}
	\draw (0,0) -- (5,0);
	\draw (0.5,0) node[label=center:{\large $\wedge$}] {};
	\draw (1,0) circle(3pt);
	\draw[dashed] (0.5,0.1) .. controls (0.625, 0.4) and (0.875,0.4) .. (1,0);
	\draw (1.5,0) node[label=center:{\large $\times$}] {};
	\draw (2,0)  circle(3pt);
	\draw (1.5,0) .. controls (1.625,0.4) and (1.875,0.4) .. (2,0);
	\draw (2.5,0) node[label=center:{\large $\times$}] {};
	\draw (3,0)  circle(3pt);
	\draw (2.5,0) .. controls (2.625,0.4) and (2.875,0.4) .. (3,0);
	\draw (3.5,0)  circle(3pt);
	\draw (4,0) node[label=center:{\large $\times$}] {};
	\draw (4.5,0)  circle(3pt);
	\draw (4,0) .. controls (4.125,0.4) and (4.375,0.4) .. (4.5,0);
\end{tikzpicture}
\]
As a $\fg_x\times\mathbb{C}\langle h\rangle$-module we have that
\[
DS_xL(\lambda)\cong L(\lambda_1)_1\oplus \Pi L(\lambda_1)_{-1}\oplus L(\lambda_2)_1\oplus \Pi L(\lambda_2)_{-1}\oplus L(\lambda_3)_2\oplus \Pi L(\lambda_3)_{-2},
\]
where $\lambda_1=7\epsilon_1+4\epsilon_2-4\epsilon_4-7\epsilon_5$, $\lambda_2=7\epsilon_1+2\epsilon_2-2\epsilon_4-7\epsilon_5$, and $\lambda_3=4\epsilon_2+2\epsilon_3-2\epsilon_5-4\epsilon_5$.  Here for a $\fg_x$-module $V$ and $t\in\mathbb{C}$ we write $V_t$ for the $\fg_x\times\mathbb{C}\langle h\rangle$-module on which $h$ acts by $t$.

\subsection{Depth}\label{depth}
We  set $X(\fg)_r:=\{x\in \fg_{\ol{1}}^{ss}|\ \rank x=r\}$.
For a $\fg$-module $N$ we set 
$$\breve{X}(N):=\{x\in \fg_{\ol{1}}^{ss}\setminus\{0\}|\ \DS_x(N)\not=0\}$$

and  introduce $\depth(N)\in\frac{1}{2}\mathbb{N}\cup\{\infty\}$ recursively by
$$ \depth(N):=\left\{\begin{array}{ll}
0 & \text{ if }\breve{X}(N)=\emptyset\\
\displaystyle\max_{x\in \breve{X}(N)}\bigl(\depth(\DS_x(N))+\rank x\bigr) & \text{ if }\breve{X}(N)\not=\emptyset.
\end{array}\right.$$

By~\Cor{corchi} one has $\depth(N)\leq r$ if  $N$ has a central character of atypicality $r$.

\subsubsection{}
\begin{cor}{cordepth}
For  a finite-dimensional simple module $L=L(\lambda)$  one has 
$$\depth(L)=\atyp \lambda.$$
\end{cor}
\begin{proof}
By~\Cor{corchi} one has $\depth(L)\leq \atyp\lambda$. Assume that $\atyp\lambda>0$.
Combining~\Thm{thmDSqs} and~\Cor{DS semisimple} we deduce
the existence of $x\not=0$ (with $\rank x\leq 1$) such that
$\DS_x(L)$ has a simple direct summand $L(\nu)$.
By~\Cor{corchi}, $\atyp\nu=\atyp\lambda-\rank x$. Using the induction
on $\atyp\lambda$ we obtain $\depth (L)\geq \atyp \lambda$.
\end{proof}

\subsection{Proof of~\Thm{thmDSqs} (ii)}\label{subsec smult}
Let us show that 
the following are equivalent
\begin{itemize}
	\item[(a)] 
	$\operatorname{smult}(\DS_xL(\lambda),L(\nu))\not=0$;
	\item[(b)] 
	$[\DS_x(L(\lambda)):L(\nu)]=1$;
	
	\item[(c)] $\zero(\lambda)-\zero(\nu)=2s$
	
\end{itemize}
where $\operatorname{smult}$ (super multiplicity)  of $L$ in $M$ is given by:
\[
\operatorname{smult}(M,L):=\left\{\begin{array}{ll}
	[M:L]-[M:\Pi L] & \text{ if }L\not\cong\Pi L\\
	{[}M:L]\text{ mod }2 & \text{ if }L\cong\Pi L\\
\end{array}\right.
\]

The implication (c) $\Longrightarrow $ (b) can be easily established; for a proof in the case of $x^2=0$, see for example~\cite{Gcore}, Prop. 5.7.2. It remains to verify that
the super multiplicity of $L(\nu)$ in $\DS_x\bigl(L(\lambda)\bigr)$ 
is zero except for the case when $\zero(\lambda)-\zero(\nu)=2s$.



\subsubsection{}\begin{cor}{ds_1/2}
	For $x$ of rank $r$ we have
	\[
	[\Res_{\fq_{n-2r}}^{\fq_n} L(\lambda)]=[\DS_x\bigl(L(\lambda)\bigr)]=\left\{\begin{array}{ll}
		[L(\lambda')] & \text{ if }\zero(\lambda)\geq 2r\\
		0 & \text{ otherwise}\\
	\end{array}\right.
	\]
	Here $\lambda'$ is the dominant weight obtained from $\lambda$ by removing $2r$ zeroes.
\end{cor}
\begin{proof}
	Take $x=C_r$ and use~\Thm{thmDSqs} (iii).
\end{proof}

\subsubsection{}\begin{cor}{}
	Let $x$ be of rank $r$, and let $\lambda\in P^+(\fq_n)$, $\mu\in P^+(\fq_{n-2r})$. 
	One has $\operatorname{smult}(\DS_x \bigl(L(\lambda)\bigr),L(\mu))=\pm 1$	
	if $\mu$ is  obtained from $\lambda$ by removing $2r$ zeroes and
	$\operatorname{smult}(\DS_x \bigl(L(\lambda)\bigr),L(\mu))=0$	otherwise.
\end{cor}
\begin{proof}
In $\Gr_-(\fq_n)$ we have the following identity:
\[
[M]=\sum\limits_{L}\operatorname{smult}(M,L)[L],
\]
where $L$ runs over irreducible representations of $\fq_n$ up to parity. 
Now the statement follows from ~\Cor{ds_1/2}.\end{proof}

This completes the proof~\Thm{thmDSqs} (ii).

\subsection{Examples}\label{example comps}
\begin{enumerate}
	\item For $\lambda:=(4,1,0,-1,-4)$ we obtain the arc diagram: 
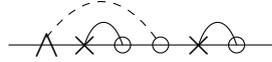
\[
\begin{tikzpicture}
	\draw (0,0) -- (3.5,0);
	\draw (0.5,0) node[label=center:{\large $\wedge$}] {};
	\draw (1,0) node[label=center:{\large $\times$}] {};
	\draw (1.5,0)  circle(3pt);
	\draw (1,0) .. controls (1.125,0.4) and (1.375,0.4) .. (1.5,0);
	\draw (2,0)  circle(3pt);
	\draw[dashed] (0.5,0.1) .. controls (1, 0.75) and (1.5,0.75) .. (2,0);
	\draw (2.5,0) node[label=center:{\large $\times$}] {};
	\draw (3,0)  circle(3pt);
	\draw (2.5,0) .. controls (2.625,0.4) and (2.875,0.4) .. (3,0);
\end{tikzpicture}
\]
We see that the diagram $\Arc(\lambda)$ has three arcs 
$$\arc(0;3), \  \arc(4;5), \ \arc(1;2);$$
and that the first two arcs are maximal, and the first is a half arc.  If $x$ is of rank $1/2$, then  $\DS_{x_1}(L(\lambda))=L(\lambda')$ where $\lambda'=(4,1,-1,-4)$.  If $x$ is of rank $1$ with $x=T_{0,B}$ and $\tr(B)=0$, then $\DS_x(L(\lambda))$ has length two with subquotients $L(\nu), \Pi L(\nu)$ for $\nu=(1,0,-1)$. If $x$ is of rank $1$ with $x=T_{0,B}$ and $\tr(B)\neq 0$, then $\DS_x (L(\lambda))=0$.  

\item By~\Cor{cora}, for an atypical weight
$\lambda\in P^+(\fq_3)$  one has $\DS_x(L(\lambda))=0$ for any $x$ with $x^2=0$ if and only if 
$\lambda=(1,0,-1)$; for  $\fq_4$ this holds  if and only if $\lambda$ 
$$
\lambda\in\{ (2,0,-1,2),\ (2,1,0,-2),\ (a,1,0,-1),\ (1,0,-1,-a)\}\ \text{ for }a\in\mathbb{N}_{>1}.$$

\item
A weight $\lambda\in  P^+(\fq_5)$  of atypicality $2$ satisfies $\DS_x(L(\lambda))=0$ for any $x$ with $x^2=0$ if and only if 
$\lambda=(2,1,0,-1,-2)$ or $\lambda=(3,1,0,-1,-3)$.
In these cases $\DS_{1/2}(L(\lambda))$ is  $L(\nu)$ for $\nu=(2,1,-1,-2)$ and $\lambda=(3,1,-1,-3)$ respectively.

\item Let $\lambda$ be one of the weights $(2,1,0,-1,-2),(3,1,0,-1,-3)$.  By Lemma \ref{ds dim bound lemma}, if $x=x_{ss}+x_{nil}\in\fg_{\ol{1}}^{ss}$ satisfies $\DS_xL(\lambda)\neq0$, then $x_{nil}=0$.  

\end{enumerate}

\subsection{Shrinking}\label{shrinking}
Let $f$ be a core-free diagram. Each minimal arc in $\Arc(f)$ takes the following form $\arc(a;a+1)$ for $a\geq0$ such that $f(a)=\times$ or $\wedge$, and $f(a+1)=\circ$.

If $a$ supports a minimal arc
 we define the diagram  $\shr_a(f)$ as follows:

\begin{itemize}
\item
if  $a\not=0$,  we produce
$shr_a(f)$  by ``shrinking''  the positions $a$ and  $a+1$;
for instance,
$$\shr_3(\wedge^2\circ\times\times\circ\times)=\wedge^2\circ\times\times;$$

\item if $\arc(0;1)\in\Arc(f)$ we produce
$shr_a(f)$  by a removal of a 
$\wedge$ from the zero position and ``shrinking'' the position
$1$; for instance, 
\[
\shr_0(\wedge^3\circ\times\circ\times)=\wedge^2\times\circ\times.
\]
\end{itemize}

Observe that 
$\Arc(shr_a(f))$ is obtained from $\Arc(f)$ by ``shrinking'' the minimal arc
supported at $a$; the rest of the arcs remain ``the same'', i.e.
we have the natural injective map $\Arc(shr_a(f))\to \Arc(f)$
which preserves order between arcs.

The following proposition, 
which will be proven in Section~\ref{verblud}, 
reduces the computation of $\DS_x(L(\lambda))$ to
the case $\fg=\fq_{s}$, where $\rank x=s$.

\subsection{}
\begin{prop}{corshrink} Let $x\in\fg_{\ol{1}}^{ss}$ be of rank $s$, and let $\lambda\in P^+(\fq_n)$, $\nu\in P^+(\fq_{n-2s})$.
\begin{enumerate}
\item
If $m_x(\lambda;\nu)\not=0$, then $\core(\lambda)=\core(\nu)$ and
 $\Arc(\nu)\subset \Arc(\lambda)$. 
\item
If $\core(\lambda)=\core(\nu)$ and
 $\Arc(\nu)\subset \Arc(\lambda)$, then
\begin{equation}\label{DShowl}
	m_x(\lambda;\nu)=\dim \DS_x(L(\lambda'))),
\end{equation}
where 
$\Arc(\lambda')$ is obtained from $\Arc(\howl(\lambda))$
by shrinking successively all arcs appearing in $\Arc(\howl(\nu))$ (starting from the minimal arcs).\end{enumerate}
\end{prop}

\subsubsection{Example} Suppose that $\rk x=3$.  Then \Prop{corshrink} gives
\[
[\DS_x(L(\wedge^4\circ\times\times\circ\times):L(\wedge^2\circ\circ\times)]=0
\] 
since $\arc(0,2)$ lies in $\Arc(\wedge^2\circ\circ\times)$ but does not lie in
$\Arc(\wedge^4\circ\times\times\circ\times)$.
Another example: for $\rk x=2$ we have
\[
[\DS_x(L(\wedge^2>\circ\times\circ\circ\times\times):L(\wedge^2>\circ\times)]=\dim\DS_x(L(\circ\times\times)).
\]
The above formula is obtained as follows: in this case
$\howl{\nu}$ has the diagram $f=\wedge^2\circ\times$
with minimal arcs $\arc(0,1)$ and $\arc(2;3)$; shrinking these arcs gives the diagram
$\shr_0\shr_2(f)=\wedge$, which has  a unique  $\arc(0;1)$.
The diagram of $\howl{\lambda}$ is $g=\wedge^2\circ\times\circ\circ\times\times$
and 
\[
\shr_0\shr_0\shr_2(g)=\shr_0\shr_0(\wedge^2\circ\circ\times\times)=\circ\times\times.
\]

\subsection{Reduction of Theorems~\ref{thmDSqs}--\ref{thm rk 1 half} and~\Thm{thm grading} to~\Prop{corshrink} ??}
\Prop{corshrink} implies~\Thm{thmDSqs} (i) and (ii) is proven in Section \ref{subsec smult}.
In Theorems~\ref{main thm rk 1}, \ref{thm rk 1 half}, $x$ is of rank $1$.
For this case  the weight
  $\lambda'$ appearing in the right-hand side of~(\ref{DShowl})
  is a core-free weight
of $\fq_{2}$, so~\Prop{corshrink} (ii)  reduces the assertions to the case 
$\fg=\fq_2$, which can be easily computed explicitly:
we have to compute $\DS_x(L(\lambda'))$ for a weight $\lambda'\in P^+(\fq_2)$ of atypicality 2 and $x\in (\fq_2)_{\ol{1}}^{ss}$ of rank $1$.
 Recall that $\DS_x(L(\lambda'))=0$ if $\lambda'$ is not integral and $x^2\not=0$.
 One has
\[
\dim\DS_x(L(\lambda'))=\left\{\begin{array}{lll}
1 & \text{ if }\lambda'=0 & \text{ i.e. } \diag\lambda'=\wedge^2\\
2 & \text{ if }\lambda'\not=0 & \text{ i.e. }\diag\lambda'=\circ\ldots\times\\
\end{array}\right.
\]
if  $x^2=0$ 
  or if  $\lambda'$ is integral and
   $x=T_{0,B}$ with $\tr(B)=0$. This establishes~Theorems~\ref{main thm rk 1} (i) and \ref{thm rk 1 half}.  In the remaining case when
   $\lambda'$ is integral and 
$x=T_{0,B}$ with $B$ of rank 2 with $B^2$ semisimple and
 $\tr(B)\neq 0$ we have
\[
\dim\DS_x(L(\lambda'))=\left\{\begin{array}{lll}
	1 & \text{ if }\lambda'=0 & \text{ i.e. } \diag\lambda'=\wedge^2\\
	0 & \text{ if }\lambda'\not=0 & \text{ i.e. }\diag\lambda'=\circ\ldots\times\\
\end{array}\right.
\]
This establishes~\Thm{main thm rk 1} (ii).
Finally, (iii) of~\Thm{thmDSqs} follows from the following lemma:

\subsubsection{}\begin{lem}{ss induction lemma} 
	We have $\DS_{C_n}(L(\lambda))=0$ for a non-zero weight $\lambda\in P^+(\fq_n)$.
\end{lem}

\begin{proof} 
Recall that $\DS_x(N)$ is an $\fg^x$-subquotient $\Res_{\fg^x}^{\fg} N$
which is annihilated by the ideal $[x,\fg]\cap \fg^x$. For $x:=C_n$ one has $\fg^x=[x,\fg]=\fg_{\ol{0}}$, so 
  $\DS_{C_n}(N)$  is a subquotient of $N^{\fg_{\ol{0}}}$ if
  $N$ is a finite-dimensional $\fg$-module. 
By~\cite{Cheng}, $\Res_{\fg_{\ol{0}}}^{\fg} L(\lambda)$ contains no copies of the trivial module if $L$ is not itself trivial; thus we must have $\DS_{C_n}(L(\lambda))=0$.
\end{proof}

\subsubsection{Reduction of Theorem \Thm{thm grading}}\label{section thm grading summary} In order to prove Theorem \ref{thm grading}, we will follow the proof of Proposition \ref{corshrink} given in the next section while also keeping track of the action of $h$.  Namely, suppose that we are in the setup of Proposition \ref{corshrink} except that we assume that $x$ is of rank 1 with $x^2=0$.  Then $h$ will act on any composition factor $L(\nu)$ inside of $\DS_xL(\lambda)$ by a scalar; thus we may upgrade $m_x(\lambda,\nu)$ to a vector space with an $h$-action that we write as $M_x(\lambda,\nu)$, where we have $\dim M_x(\lambda,\nu)=m_x(\lambda,\nu)$.  Then in this case we will show that we obtain an isomorphism of $h$-vector spaces
\begin{eqnarray}\label{eqn grading}
M_x(\lambda,\nu)\cong \DS_xL(\lambda').
\end{eqnarray}
The isomorphism (\ref{eqn grading}) reduces the result to the $\fq_2$-case, where we find that, for $k\in\frac{1}{2}\mathbb{N}$, (up to parity)
\[
\DS_xL(k(\epsilon_1-\epsilon_2))=\mathbb{C}_k\oplus\Pi\mathbb{C}_{-k},
\]
where $\mathbb{C}_k$ is the one-dimensional even vector space on which $h$ acts by $k\Id$.

\section{Proof of ~\Prop{corshrink}}\label{verblud}

\subsection{Translation functors and $\DS$}\label{TransDS}
Our main tools  are translation functors described in~\cite{Br}.
Below we briefly 
recall the connection between the translation functor and $\DS$ functor,
see~\cite{GH}, 7.1.

\subsubsection{Notation}
For  a core diagram $c$ we denote by $\Irr(\fg)^c$
the set of isomorphism classes of finite-dimensional irreducible modules $L(\lambda)$
with $\core(\lambda)=c$. Recall that all modules in $\Irr(\fg)^c$
have the same central character $\chi$.  Let 
$\Fin(\fg)^{c}$ be
the full subcategory  of $\Fin(\fg)$ 
which corresponds to the central character $\chi$
(i.e., $N\in\Fin(\fg)^{c}$ if and only if $(z-\chi(z))^{\dim N} N=0$
for each $z\in\cZ(\fg)$). Then $\Irr(\fg)^c$
is the set  of isomorphism classes of irreducible modules in $\Fin(\fg)^c$.

Note that $L(\vareps_1)$ is the standard module and $L(-\vareps_n)$
is its dual.
Let $c,c'$ be core diagrams; we denote by $T_{c}^{c'}$ (resp.,
$(T_{c}^{c'})^*$) 
the translation functor $\Fin(\fg)^{c}\to \Fin(\fg)^{c'}$
which maps $N$ to the projection of 
$N\otimes L(\vareps_1)$ (resp., of $N\otimes L(-\vareps_n)$)
to the subcategory $\Fin(\fg)^{c'}$.  The functors
$T_{c}^{c'}$,
$(T_{c}^{c'})^*$ are both left and right adjoint to each other.

\subsubsection{}
It is easy to check that for any $x\in\fg_{\ol{1}}^{ss}$ we have $\DS_x(L_{\fg}(\vareps_1))=L_{\fg_x}(\vareps_1)$.
Since $\DS$ commutes with tensor product and preserves cores 
(see~\ref{corchi}) one has 
$$\DS_x(T^{c'}_{c}(N))\cong T^{c'}_{c}(\DS_x(N))$$
where $c',c$ are core diagrams and $T^{c'}_{c}$ stands for the functors $\Fin(\fg)^{c}\to \Fin(\fg)^{c'}$
and $\Fin(\fg_x)^{c}\to \Fin(\fg_x)^{c'}$.

Since the translation functors are exact, they induce morphisms of the Grothendieck ring.
For any  $N\in\Fin(\fg)^{c}$ and $L'\in \Irr(\fg_x)^{c'}$
 we have
\begin{equation}\label{Tata}\begin{array}{ll}
[\DS_x(T^{c'}_{c}(N)):L']&=[T^{c'}_{c} (\DS_x(N)):L']\\
&=\displaystyle\sum_{L_1\in\Irr(\fg_x)^{c}}
[\DS_x(N):L_1][T^{c'}_{c}(L_1):L'].\end{array}\end{equation}

Assume that $L_1\in\Irr(\fg_x)^{c}$ and $L'\in \Irr(\fg_x)^{c'}$ are such that
for each $L_2\in\Irr(\fg_x)^{c}$ with $L_2\not\cong L_1,\Pi L_1$ one has
 $[T^{c'}_{c}(L_2):L']=0$. Then~(\ref{Tata}) gives 
\begin{equation}\label{Tfu}
[\DS_x(T^{c'}_{c}(N)):L']=[T^{c'}_{c}(L_1):L']\cdot [\DS_x(N):L_1].
\end{equation}

\subsection{Outline of the proof of~\Prop{corshrink}}
\label{planverbluda}
We will use the following notation
\[
\begin{array}{c}  \frac{\diag(\lambda)}{\diag(\nu)}:=[\DS_x(L(\lambda)):L(\nu)].
\end{array}
\]
By~\ref{corchi}, $\frac{f}{g}\not=0$ implies that
 $\core(f)=\core(g)$. The proof of the remaining assertions of~\Prop{corshrink} 
follows the same steps as in $\osp$-case.
We briefly describe these steps below.

\subsubsection{Reduction to the stable case}\label{stabilize}
We call a weight diagram $f$  {\em stable} if all symbols
$\times$ precede all core symbols, (notice that the symbols $\wedge$ necessarily precede core symbols); we say that $\eta\in P^+(\fg)$
is stable if $\eta$ is typical or $\diag(\eta)$ is stable. For instance,
$(4,1,0,-1,-2)$ with the diagram $\wedge\times<\circ>$
is stable and the weight $\zeta=(2,1,-2)$ 
with the diagram  $\circ>\times$ is not stable.

The general case can be reduced to 
the stable case with the help of translation functors 
described in~\cite{Br}. For this step we will use  translation 
functors which preserves $\howl$:
such functor transforms $L(\eta)$ 
to $L(\eta')\oplus \Pi L(\eta')$, where  $\diag(\eta')$ is obtained from 
$\diag (\eta)$ by permuting two neighboring symbols at non-zero positions
if exactly one of these symbols is a core symbol:
for instance, from $\zeta$ as above we can obtain $\zeta'$ with the diagram
 $\circ\times>$ which is stable.
Note that $\howl({\eta})=\howl({\eta}')$.
Using these  functors we can transform any two simple modules
$L(\lambda)$, $L(\nu)$ with $\core(\lambda)=\core(\nu)$
to the modules $L(\lambda_{st})^{\oplus r}\oplus \Pi L(\lambda_{st})^{\oplus r}$, $L(\nu_{st})^{\oplus r}\oplus \Pi L(\nu_{st})^{\oplus r}$,
where $\lambda_{st},\nu_{st}$ are stable weights with $$\core(\lambda_{st})=\core(\nu_{st}),\ \ \
\howl(\lambda)=\howl(\lambda_{st}),\ \ \ \howl(\nu)=\howl(\nu_{st})$$
(the diagrams of $\lambda_{st}$ and $\nu_{st}$
 are stable diagrams obtained from $\diag(\lambda)$, $\diag(\nu)$ by 
moving all core symbols from the  non-zero position ``far enough'' to the right).
 For instance,
$$\begin{array}{lccl}
\lambda=(6,5,1,0,0,-3,-5) & & & \nu=(6,1,0,0,-3)\\
\diag(\lambda)=\wedge^2>\circ<\circ\times> & & &
\diag(\nu)=\wedge^2>\circ<\circ\circ>\\
\diag(\lambda_{st})=\wedge^2\circ\circ\times\circ\circ><> & & &
\diag(\nu_{st})= \wedge^2\circ\circ\circ\circ\circ><>   \\
\lambda_{st}=(8,6,3,0,0,-3,-7) & & & \nu_{st}=(8,6,0,0,-7).
\end{array}$$
Using~(\ref{Tfu})  we will obtain the formula~(\ref{formulastab}) which implies
$$[\DS_x(L(\lambda)):L(\nu)]=[\DS_x(L(\lambda_{st})):L(\nu_{st})].$$
This formula reduces  the general case to the case when $\lambda,\nu$ are stable.

\subsubsection{Reduction to the core-free case}\label{step1}
In~\ref{howl} we will show that  for  a stable weight $\lambda$ the computation of multiplicities in
$\DS_x(L(\lambda))$ can be  reduced to the case when
$\lambda$ is core-free.

\subsubsection{Shrinking in the core-free case}\label{shrink}
 Using translation functors 
which do not preserve $\howl$
we will obtain useful ``cancellation'' formulae~(\ref{longformula}). Applying these formulae 
to core-free weights $\lambda,\nu$ with the diagrams $f,g$ we obtain
\begin{enumerate}
\item if 
$[\DS_x(L(\lambda)):L(\nu)]\not=0$, then each  minimal 
arc in $\Arc(\nu)$ is  a minimal arc in $\Arc(\lambda)$;
in particular, if the operation $shr_a$ (introduced in~\ref{shrinking}) is defined for $g$, then
this operation is defined for $f$;
\item
if $shr_a$ is defined for  $f$ and $g$ for $a\not=0$, then
$
\frac{f}{g}=\frac{shr_a(f)}{shr_a(g)}$;
\item
if 
$[\DS_x(L(\lambda)):L(0)]\not=0$, then $\Arc(0)\subset\Arc(\lambda)$
and $\frac{f}{g}=\frac{shr_0(f)}{shr_0(g)}$.
\end{enumerate}

\Prop{corshrink} follows from these assertions by induction on $\nu_1$
(where $\nu=(\nu_1,\ldots,\nu_{n-2r})$).

\subsection{Useful translation functors}\label{usefultr}
Below we recall 
several results of Lemma 4.3.8 in~\cite{Br}. 
Let $g$ be a weight diagram and
$c$ be the core diagram of $g$. We retain notation of~\ref{TransDS}.

\subsubsection{Translation functors preserving $\howl$}\label{transeq}
Fix $a\not=0$. 
For a weight diagram $f$ we denote by $\sigma(f)$ 
the diagram obtained by interchanging the symbols at the positions $a$ and $a+1$.

If  $c(a)=>$ and $c(a+1)=\circ$ (i.e. 
 $c=\ldots>\circ\ldots$, $\ \sigma(c)=\ldots\circ>\ldots$), then
$$T_{c}^{\sigma(c)}(L(g))=L(\sigma(g))\oplus \Pi L(\sigma(g)).$$

The same formula holds if $c(a)=\circ$, $c(a+1)=<$.
We depict these formulae as
$$\begin{array}{l}
T_{\ldots>\circ\ldots}^{\ldots\circ>\ldots}\bigl(L(\ldots> \star\ldots)\bigr)=
L(\ldots \star > \ldots)^{\oplus 2}\\
T^{\ldots>\circ\ldots}_{\ldots\circ>\ldots}\bigl(L(\ldots\star <\ldots)\bigr)=
L(\ldots <\star \ldots)^{\oplus 2}
\end{array}$$
for  $\star\in \{\circ,\times\}$. The similar formulae hold for $T^*$ if  we interchange $>$ and $<$.

\subsubsection{Translation functors reducing atypicality}\label{transred}
Using the above notation  we have 
$$\begin{array}{l}
T_{\ldots\circ\circ\ldots}^{\ldots<>\ldots}\bigl(L(\ldots\times \times\ldots)\bigr)=
T_{\ldots\circ\circ\ldots}^{\ldots<>\ldots}\bigl(L(\ldots\circ \star\ldots)\bigr)=0
\\
T_{\ldots\circ\circ\ldots}^{\ldots<>\ldots}\bigl(L(\ldots\times \circ\ldots)\bigr)=L(\ldots < > \ldots)^{\oplus 2}\\
T_{\circ\circ\ldots}^{\circ>\ldots}\bigl(L(\wedge^r \times\ldots)\bigr)=T_{\circ\circ\ldots}^{\circ>\ldots}\bigl(L(\circ \star\ldots)\bigr)=0
\\
T_{\circ\circ\ldots}^{\circ>\ldots}\bigl(L(\wedge^{r+1}\circ\ldots)\bigr)=
L(\wedge^r>\ldots)^{\oplus d},
\end{array}
$$
where  $\star\in \{\circ,\times\}$,  $r\geq 0$ and $d=1$ if $n-r$ is odd, $d=2$ if $n-r$ is even. For instance, 
$$T_{\circ\hspace{0.27em} \circ\hspace{0.27em}  \circ \hspace{0.1em} >>}^{\circ\hspace{0.1em} < >>> }(L(4,3,1,-1))=L(4,3,2,-1)^{\oplus 2},$$ 
since $(4,3,1,-1)$  has the diagram  $\circ\times\circ>>$ 
$(4,3,2,-1)$ has the diagram $\circ<>>>$.

\subsection{Cancellation formulae}
We denote by $f_-f_+$
the diagram obtained by "gluing"
the diagrams $f_-$ and $f_+$ (where $f_+$ has exactly one symbol at each position)
for instance,
$$\begin{array}{l}
f_-=\wedge^4\circ\ \ \  f_+=\times,\ \ \ \ f_-f_+=\wedge^4\circ\times.
\end{array}
$$

Recall that $\frac{g}{g'}$ stands for the multiplicity
$[\DS_s(L(\lambda)):L(\nu)]$, where $\diag\lambda=g$ and $\diag\nu=g'$.
Applying~(\ref{Tfu}) to the translation functors which appeared in~\ref{transeq},  we obtain the following formulae
\begin{equation}\label{formulastab}
\begin{array}{llr}
\frac{f_->\star f_+}{g_->\star g_+}=\frac{f_-\star >f_+}{g_-\star >g_+} &
 & \frac{f_-<\star f_+}{g_-<\star g_+}=\frac{f_-\star <f_+}{g_-\star <g_+},
\end{array}
\end{equation}
where the symbol $\star\in\{\circ,\times\} $ occupies the same position in each diagram
(for instance, the first formula gives
 $\frac{\circ>\times\times\times}{\times>\times\circ\circ}=\frac{\circ\times>\times\times}{\times\times>\circ\circ}$).
Using~(\ref{Tfu}) for the translation functors appearing
in~\ref{transred} we get for $r,r'\geq 0$

\begin{equation}\label{longformula}
\begin{array}{llr}
\frac{f_-\circ \star f_+}{g_-\times\circ g_+}=\frac{f_-\times\times f_+}{g_-\times\circ g_+}=0 & &
 \frac{\wedge^r\times f_+}{\wedge^{r'+1}\circ
g_+}=0 \ \ \\
\frac{f_-\times\circ f_+}{g_-\times\circ g_+}=\frac{f_-<>f_+}{g_-<>g_+} & &
\frac{\wedge^{r+1}\circ f_+}{\wedge^{r'+1}\circ g_+}=
\frac{\wedge^{r}>f_+}{\wedge^{r'}>g_+}.
\end{array}
\end{equation}

\subsection{Reduction to the core-free case}\label{howl}
Our next goal is to verify 
$$\frac{f}{g}=\frac{\howl(f)}{\howl(g)}$$
if $f,g$ are stable diagrams with $\core(f)=\core(g)=c$.
The formula is tautological if $f,g$ are core-free, so we assume that
this is not the case, i.e. $c$ is non-empty.

\subsubsection{Notation for~\ref{howl}}
Let $\lambda\in \ft^*$, $\nu\in\ft_x^*$ be stable weights with 
$$\core(\lambda)=\core(\nu)=c.$$

Let $u_+$ (resp., $u_-$) be the number of symbols $>$
(resp., $<$) of $c$; we set
$$u:=u_-+u_+,\  \ \  I':=\{u_++1,\ldots, n-u_-\},\ \ I'':=I\setminus I'.$$
The assumption on $c$  implies $u\not=0$.
Since $\DS$ ``preserves the cores'', for $x\in\fg_{\ol{1}}^{ss}$ of rank $s$ we have $\DS_x(\Fin(\fg)^c)=0$
if $n-2s<u$. Therefore we assume
$$0<u\leq n-2s.$$

Retain notation of~\ref{notatqJ} and introduce
$$\fg':=\fq(I'),\ \ \fh':=\fh(I'), \ \ \fh'':=\fh(I''),\ \ \ft':=\ft(I'), \ \ \ft'':= \ft(I'').$$
Note that $\fg'\cong \fq_{n-u}$, $\fh'$ is a Cartan subalgebra of $\fg'$ and  $\fh=\fh'\times\fh''$.

Consider the subalgebra
$\fg'+\fh=\fg'\times \fh''\subset\fg$.
We will use the notation  $L_{\fg'\times\fh''}(\mu)$, $L_{\fg'}(\mu')$
for the corresponding simple modules over $\fg'\times\fh''$  and $\fg'$ respectively
(where  $\mu\in\ft^*$  and $\mu'\in (\ft')^*$).

Let $\wt(c)\in (\ft'')^*$ be the weight ``corresponding to the diagram $c$'':
we take $\langle\wt(c),h_i\rangle$ (resp., $\langle\wt(f),h_{n+1-i}\rangle$) 
 equal to the coordinate of $i$th symbol $>$ (resp., $<$) in $c$
counted from the right. 
For example, for $c=>\circ\circ><>$ and $n=8$, $\ft'$ (resp.,
$\ft''$) is spanned
by $h_3,\ldots,h_7$ (resp., 
$h_1,h_2, h_{8}$) and $\wt(c)=5\vareps_1+3\vareps_2-4\vareps_8$.

\subsubsection{Remark}\label{remstable}
For $\mu\in P^+(\fq_n)$ with $\core(\mu)=c$ the following are 
equivalent
\begin{itemize}
\item $\mu$ is stable;
\item
$\mu|_{\ft'}=\howl(\mu)$;
\item $\mu|_{\ft''}=\wt(c)$.
\end{itemize}

\subsubsection{Choice of $x$}
We take 
$x$ of rank $s$ lying in $\fg'$ (this can be done since $2s\leq n-u$). One has
$\fg'_x:=\DS_x(\fg')\cong \fq_{n-u-2s}$. 
We identify $\fg_x$ (resp., $\fg'_x$) with the subalgebras 
of $\fg=\fq_n$ as in~\ref{DSqident} and
identify $\DS_x(\fg'\times \fh'')$ with $\fg'_x\times\fh''\subset \fg_x$.

\subsubsection{}
\begin{lem}{lemele}
$[\DS_x(L(\lambda)):L_{\fg_x}(\nu)]=
[\DS_x(L_{\fg'\times\fh''}(\lambda)):L_{\fg'_x\times\fh''}(\nu)]$.
\end{lem}
\begin{proof}
Choose $z\in\ft''\subset\ft$ in such a way  that $\{\vareps_i(z)\}_{i\in I''}$
are positive real numbers linearly independent over $\mathbb{Q}$. We set
$$a_0:=\wt(c)(z).$$
One has
\begin{equation}
\label{gz}
\fg^z=\fg'\times \fh'',\ \ \ \ \lambda(z)=\nu(z)=a_0
\end{equation}
(the last formula follows from~\ref{remstable} and the stability of $\lambda,\nu$).

For a $\ft''$-module $N$ we denote by $\Spec_z(N)$ the set of $z$-eigenvalues on
$N$,  by $N_a$  the $a$th eigenspace and view $N_a$ as a module over $\fg^z$.

It is easy to see that for any $\mu\in\ft^*$ we have
\begin{itemize}
\item[(a)]
 $\Spec_z(L(\mu))\subset \mu(z)-\mathbb{R}_{\geq 0}$;
\item[(b)] $L(\mu)_{\mu(z)}=L_{\fg'\times\fh''}(\mu)\ $ (this  follows from the PBW-theorem).
 \end{itemize}

Set
$N:=\DS_x(L(\lambda))$. 
Since $x\in\fg'$ we have $[z,x]=0$ so $N_{a}=\DS_x\bigl( L(\lambda)_{a}\bigr)$; this gives
$$\Spec_z\bigl(N)\subset  \lambda(z)-\mathbb{R}_{\geq 0},\ \ \ \ N_{a_0}=\DS_x\bigl(L_{\fg'\times\fh''}(\lambda) \bigr).$$
In particular,
\begin{equation}\label{NLNL}
[\Res^{\fg_x}_{\fg_x'\times\fh''} N:L_{\fg'_x\times\fh''}(\nu)]=[N_{a_0}:
L_{\fg'_x\times\fh''}(\nu)]=[\DS_x\bigl(L_{\fg'\times\fh''}(\lambda) \bigr):
L_{\fg'_x\times\fh''}(\nu)].\end{equation}
On the other hand,
\begin{equation}\label{NLNL2}
[\Res^{\fg_x}_{\fg_x'\times\fh''} N:L_{\fg'_x\times\fh''}(\nu)]=
\sum_{\mu\in \ft_x^*} [N:L_{\fg_x}(\mu)]\cdot [L_{\fg_x}(\mu) :L_{\fg'_x\times\fh''}(\nu)].\end{equation}

Assume that 
$[N:L_{\fg_x}(\mu)]\cdot [L_{\fg_x}(\mu) :L_{\fg'_x\times\fh''}(\nu)]\not=0$.
Then, by above, $\nu(z)\leq \mu(z)\leq a_0$. Using~(\ref{gz}) we get  $\mu(z)=\nu(z)$. Applying
(b) we obtain $\mu=\nu$ and
$[L_{\fg_x}(\nu) :L_{\fg'_x\times\fh''}(\nu)]=1$. Hence~(\ref{NLNL2}) can be rewritten as 
 $$[\Res^{\fg}_{\fg'\times\fh''} N:L_{\fg'_x\times\fh''}(\nu)]=
 [N:L_{\fg_x}(\nu)] =[\DS_x(L(\lambda)):L_{\fg_x}(\nu)].$$
Now the required assertion follows from~(\ref{NLNL}).
\end{proof}

 \subsubsection{}
 \begin{lem}{lemu+12}
One has
 $$\Res^{\fg'\times\fh''}_{\fg'} L_{\fg'\times\fh''}(\lambda)=L_{\fg'}(\lambda')^{\oplus d}$$
 where $\lambda':=\lambda|_{\ft'}$ and $d:=2^{\left\lfloor\frac{u+1}{2}\right\rfloor}$.
 \end{lem}
\begin{proof}
Retain notation of~\ref{Clifford}.
Recall that $B_{\lambda}$  is the symmetric form on $\fh_{\ol{1}}$ 
given by $(H,H')\mapsto \lambda([H,H'])$. We denote by $B'_{\lambda}$
 the restriction of $B_{\lambda}$ to  $\fh'_{\ol{1}}$ and
view
$$\Cl(\lambda'):=\Cl(\fh'_{\ol{1}}, B'_{\lambda})$$
as a subalgebra of $\Cl(\lambda)$.
 By~\ref{Clifford}, the highest weight space of
$L_{\fg'\times\fh''}(\lambda)$ is  $C_{\lambda}$
and the highest weight space of
$L_{\fg'}(\lambda')$ is a simple $\Cl(\lambda')$-module which we denote by 
$E'_{\lambda'}$. It is easy to see that 
all simple subquotients of 
$\Res^{\fg'\times\fh''}_{\fg'} L_{\ft}(\lambda)$ 
are isomorphic to $L_{\fg'}(\lambda')$. 
Therefore it is enough to check that 
$$\Res^{\Cl(\lambda)}_{\Cl(\lambda')} C_{\lambda}={E'_{\lambda'}}^{\oplus d}$$
for $d$ as above.
By~\cite{GSS}, Prop. 3.5.1 the semisimplicity of $\Res^{\Cl(\lambda)}_{\Cl(\lambda')} C_{\lambda}$
is equivalent to the formula
$\Ker B'_{\lambda}=\fh'_{\ol{1}}\cap\Ker B_{\lambda}$.
This formula holds since
 $\Ker B_{\lambda}$ is spanned by $H_i$ with $\lambda_i=0$,
$\fh'_{\ol{1}}$  is spanned by $H_i$ with $i\in I'$, and
$\Ker B'_{\lambda}$ is spanned by $H_i$ with $\lambda_i=0$ and $i\in I'$. Hence
$\Res^{\Cl(\lambda)}_{\Cl(\lambda')} C_{\lambda}$ is semisimple.
By~\ref{Clifford},
$\frac{\dim C_{\lambda}}{\dim E'_{\lambda'}}=2^j$ where 
\[
j=\left\lfloor\frac{\rank B_{\lambda}+1}{2}\right\rfloor-\left\lfloor\frac{\rank B'_{\lambda}+1}{2}\right\rfloor=
\left\lfloor\frac{\nonzero(\lambda)+1}{2}\right\rfloor-\left\lfloor\frac{\nonzero(\lambda')+1}{2}\right\rfloor.
\]

By~\ref{remstable} $\ \lambda'=\howl(\lambda)$, so 
$\nonzero(\lambda')$ is even and 
$\nonzero(\lambda)=\nonzero(\lambda')+u$. This gives
$j=\left\lfloor\frac{u+1}{2}\right\rfloor$ and  completes the proof.
\end{proof}

%
%
%

\subsubsection{}
\begin{cor}{corsta}
$\frac{f}{g}=\frac{\howl(f)}{\howl(g)}$.
\end{cor}
\begin{proof}
Notice that $\ft'_x:=\ft'\cap \ft_x$
plays the same role for
$\fg_x$ as $\ft'$ for $\fg$; therefore~\ref{lemu+12} gives
$$\Res^{\fg'_x\times\fh''}_{\fg'_x} L_{\fg'_x\times\fh''}(\nu)=L_{\fg'_x}(\nu')^{\oplus d}$$
where $\nu':=\nu|_{\ft'_x}$.
Combining Lemmatta~\ref{lemu+12}, \ref{lemele} and the above formula 
we get
\[
[\DS_x(L(\lambda)):L(\nu)]=
[\DS_x(L_{\fg'}(\lambda'):L_{\fg'_x}(\nu')].
\]
By~\ref{remstable} we have  
$\lambda|_{\ft'}=\howl(\lambda)$ and $\nu|_{\ft'_x}=\howl(\nu)$
as required.
\end{proof}

\subsection{The core-free case}
It remains to verify (i)---(iii) in~\ref{shrink}.

Thus let $x\in\fg_{\ol{1}}^{ss}$ be of rank $s$, and let $\lambda\in P^+(\fq_{m+2s}), \nu\in P^+(\fq_m)$
 be core-free weights
with $m\not=0$. 

\subsubsection{}
Consider the case when $\nu\not=0$. Then $\nu_1>0$.
Write $\diag(\nu)=g_-\times\circ$ (where
 $\times$ is at the position
$\nu_1$). By~(\ref{longformula}), if
$\ [\DS_s(L(\lambda)):L(\nu)]\not=0$, then $\diag(\lambda)=f_-\times\circ f_+$
for some diagrams $f_-,f_+$ (where  $\times$ is at the position
$\nu_1$). This gives (i) in~\ref{shrink}.
Combining~(\ref{longformula}), (\ref{formulastab})  and~\ref{howl} we get
\[
\frac{f_-\times\circ f_+}{g_-\times\circ}=\frac{f_-<> f_+}{g_-<>}=\frac{f_-f_+}{g_-}.
\]
Note that both
 $\Arc(f_-\times\circ f_+)$, $\Arc(g_-\times\circ)$ contain  a minimal 
$\arc(\nu_1;\nu_1+1)$  and 
$$f_-f_+=shr_{\nu_1}(f_-\times\circ f_+),\ \ \ g_-=shr_{\nu_1}(g_-\times\circ).$$
This establishes (ii) in~\ref{shrink}.

\subsubsection{}
The same argument shows that
$[\DS_x(L(\lambda)):L_{\fq_{m}}(0)]\not=0$
implies  $\diag\lambda=\wedge^{r}\circ f_+$ for $r\geq 1$ and that
\begin{equation}\label{shr0even}
\frac{\wedge^{r}\circ f_+}{\wedge^{m}\ \ \ \ }=\frac{\wedge^{r-1}f_+}{\wedge^{m-1}\ \ \ }.
\end{equation}
Note that $\wedge^{m}=shr_0(\wedge^{m+1})$ and
$\wedge^r f_+=shr_0(\wedge^{r+1}\circ f_+)$.

From~(\ref{shr0even}) we conclude that for $m>1$ one has
\[
[\DS_x(L(\lambda)):L_{\fq_{m}}(0)]\not=0\ \ \Longrightarrow\ \ \diag\lambda=\wedge^{r}\circ f_+\ 
\] and
\[
\frac{\wedge^{r}\circ f_+}{\wedge^{m}}=\frac{\wedge^{r-1}\circ f_+}{\wedge^{m-1}}.
\]
This establishes~\ref{shrink} (iii) and completes the proof of~\Prop{corshrink}.
\qed

\subsection{Keeping track of the action of $h$}  

We now explain how to prove the isomorphism (\ref{eqn grading}) in the case when $x^2=0$ and the rank of $x$ is 1.  In this case, we view $\DS_x$ as a functor from the category of $\fg$-modules to the category of $\fg_x\times\mathbb{C}\langle h\rangle$-modules.  If $V$ is semisimple over $\fg_{\ol{0}}$, which is the only case we consider, then $\DS_x$ will admit a semisimple action of $h$.  Given a $\fg_x$-module $V$ and $t\in\mathbb{C}$, we write $V_t$ for the $\fg_x\times\mathbb{C}\langle h\rangle$-module with $h$ acting by $t$.

For the standard module $L(\epsilon_1)$, the $h$ action on $\DS_xL(\epsilon_1)$ is trivial, and thus for any pair of cores $c,c'$ we have a natural isomorphism of $\fg_x\times\mathbb{C}\langle h\rangle$-modules
\[
\DS_x(T^{c'}_{c}(N))\cong T^{c'}_{c}(\DS_x(N)).
\]
Notice that the central characters of $\fg_x\times\mathbb{C}\langle h\rangle$ are parametrized by pairs $(c,t)$, where $c$ is a core and $t\in\mathbb{C}$.  On the RHS of the above formula, $T_{c'}^{c}$ denotes the translation functor between modules with central characters of the form $(c',t)$ to modules with central characters of the form $(c,t)$.  In particular, $T_{c'}^{c}$ takes modules with central character indexed by $(c',t)$ to modules with central character indexed by $(c,t)$.

In particular we obtain the following: for any  $N\in\Fin(\fg)^{c'}$ and $L'\in \Irr(\fg_x)^{c}$
we have
\begin{equation}\label{Tata2}\begin{array}{ll}
		[\DS_x(T^{c}_{c'}(N)):L'_{t}]&=[T^{c}_{c'} (\DS_x(N)):L'_t]\\
		&=\displaystyle\sum_{L_1\in\Irr(\fg_x)^{c}}
		[\DS_x(N):(L_1)_t][T^{c}_{c'}(L_1):L'].\end{array}
\end{equation}
Using Equation (\ref{Tata2}), one may use the same operations with translation functors as is done in the general case above, and we see that the weight $t$ of $h$ will remain unaffected.

The other check that needs to be made is that for stable weights $\lambda\in P^+(\fq_n)$, $\mu\in P^+(\fq_n)$, and for $t\in\mathbb{C}$, we have
\begin{eqnarray}\label{eqn howl graded}
[DS_xL(\lambda):L(\mu)_t]=[DS_xL(\howl(\lambda)):L(\howl(\mu))_t].
\end{eqnarray}
However following the argument of Lemma \ref{lemele}, we see that the main step is to take an eigenspace of a certain semisimple operator $z$, which clearly commutes with $h$ since they both lie in $\mathfrak{t}$.  Thus Lemma \ref{lemele} becomes, in our case,
\[
[\DS_x(L(\lambda)):L_{\fg_x}(\nu)_t]=
[\DS_x(L_{\fg'\times\fh''}(\lambda)):L_{\fg'_x\times\fh''}(\nu)_t].
\]
Now using the statement of Lemma \ref{lemu+12} along with this, we obtain Equation (\ref{eqn howl graded}).  

From this, we may use the same algorithm implicitly described in Section \ref{shrink} to finally obtain the isomorphism (\ref{eqn grading}).

\section{Appendix}
The $\DS$-functor was introduced in \cite{DS}; see also~\cite{GHSS} for an expanded exposition.
We will be using a slight extension of $\DS$-functor which we define below. 

\subsection{Construction}
Let $\fg$ be a finite-dimensional Lie superalgebra.  Define
\[
\fg_{\ol{1}}^{ss}=\{x\in\fg_{\ol{1}}|\ \operatorname{ad}[x,x]\text{ is semisimple}\}.
\]

For a $\fg$-module $M$ and $x\in\fg$ we set $M^x:=\Ker_M x$.  Now let $x\in\fg_{\ol{1}}^{ss}$ and write $x^2:=\frac{1}{2}[x,x]$.  For a $\fg$-module $M$ on which $x^2$ acts semisimply we set
\[
M_x=\DS_x(M):=\frac{\ker(x:M^{x^2})}{\im(x:M^{x^2})}.
\]

Then $\DS_x$ is a tensor functor; in particular $\fg_x$ will be a Lie superalgebra, and $M_x$ will have the natural structure of a $\fg_x$-module.  Further there are canonical isomorphisms 
\[
\DS_x(\Pi(N))\cong\Pi(\DS_x(N))\ \text{ and } \ \DS_x(M)\otimes\DS_x(N)\cong\DS_x(M\otimes N).
\]
In addition we have a canonical isomorphism of $\fg_x$-modules $\DS_x(N^*)\cong (\DS_x(N))^*$.

Thus $\DS_x: M\mapsto \DS_x(M)$ is a tensor functor from the category of $\fg$-modules on which $x^2$ acts semisimply to
the category of $\fg_x$-modules.  

\subsubsection{}\label{identify}
We say that $\fg_x$  {\em can be idenitified with} a certain subalgebra of $\fg$ 
if 
\[
\fg^x=\fg_x\ltimes ([x,\fg]\cap \fg^x)
\]

\subsubsection{Grading on $\DS_x$}\label{section grading DS}  Let $x\in\fg_{\ol{1}}^{ss}$ with $x^2=0$, and write $\mathfrak{n}(x)$ for the normalizer of $x$ in $\fg$.  Then if $M$ is a $\fg$-module, $M_x$ will have the natural structure of a $\mathfrak{n}(x)$-module.  The action map gives an exact sequence
\[
0\to\fg^x\to\mathfrak{n}(x)\to\mathbb{C}\langle x\rangle,
\]
which defines a short exact sequence
 whenever there exists $h\in\fg$ such that $[h,x]=1$.

If $\fg$ is a Kac-Moody superalgebra or $\fq_n$, $\fsq_n$, such an element $h$ always exists, and further we may choose it such that under the embedding $\fg_x\subseteq\fg^x$, $h$ commutes with $\fg_x$.  In this way we obtain naturally an action of 
$\fg_x\times\mathbb{C}\langle h\rangle$ on $M_x$.  

The action of $h$ gives rise to a grading on $M_x$ as a $\fg_x$-module according to the eigenvalues of $h$:
\[
M_x=\bigoplus\limits_{t\in\mathbb{C}}(M_x)_t.
\]

\subsubsection{Hinich Lemma}
Each short exact sequence of $\fg$-modules with semisimple action of $x^2$
\[
0\to M_1\to N\to M_2\to 0
\]
induces a long exact  sequence of $\fg_x$-modules
\[
0\to Y\to \DS_x(M_1)\to\DS_x(N)\to\DS_x(M_2)\to \Pi(Y)\to 0,
\]
where $Y$ is a some $\fg_x$-module; identifying $M_1$ with its image
in $N$ we have
\[
Y=(M_1^x\cap [x,N])/(M_1^x\cap [x,M_1]).
\]
If $[x,x]=0$ and we have an element $h$ as in Section \ref{section grading DS}, the morphism $\DS_x(M_2)\to\Pi Y$ has weight $1$ for the action of $h$, while all other maps commute with $h$.

\subsubsection{$\DS$ and restriction map}\label{DSandds}
Let $\Gr(\fg)$ denote the Grothendieck ring of finite-dimensional $\fg$-modules and let $\Gr_-(\fg)$ denote the
the quotient of $\Gr(\fg)$ by the relations
$[N]=- [\Pi N]$, where $\Pi$ stands for the parity change functor.
For a finite-dimensional $\fg$-module we denote by $[N]$ its image in
$\Gr_-(\fg)$.
Although $\DS_x$ is not an exact functor, by the Hinich lemma $\DS_x$ defines  ring homomorphisms 
on reduced Grothendieck rings
$\Gr_-(\Fin(\fg))\to\Gr_-(\Fin(\fg^x))$ which coincides with the restriction map
$[M]\mapsto [\Res_{\fg^x}^{\fg}M]$, see~\cite{GHSS}.
The  ring $\Gr_-(\Fin(\fg_x))$ is a subring
of $\Gr_-(\Fin(\fg^x))$ and the image 
of $\Gr_-(\Fin(\fg))\to\Gr_-(\Fin(\fg^x))$ lies in this subring, see~\cite{GHSS}.

If $\fg_x$ is idenitified with a subalgebra of $\fg$  we have 
$$[\DS_x(M)]=[\Res_{\fg_x}^{\fg}M].$$

\subsubsection{Remark}
The morphism $[M]\mapsto[DS_xM]=[\Res_{\fg_x}^{\fg}M]$ is often denoted by $ds_x$ in many papers, including~\cite{GHSS}.  We chose to avoid this notation to emphasize the simplicity and universality of restriction, in favor of the more limited setting of the $\DS$ functor.

\subsection{Case of commuting $x,y\in \fg^{ss}_{\ol{1}}$}
Fix $x,y\in \fg^{ss}_{\ol{1}}$ and $h\in\fg$ such that 
$$[x,y]=0,\ \ [h,x ]=c_xx,\ \ \ [h,y]=c_y y\ \text{ with }\ c_x,c_y\in\mathbb{C}, \ c_x\not=c_y.$$
Note that
$x\in \fg^y$; we denote by
$\ol{x}$  the image of $x$ in $\fg_y$.

\subsubsection{}
\begin{lem}{lemdsxy}
Let $N$ be a finite-dimensional 
$\fg$-module with a diagonal action of $h,x^2,$ and $y^2$.  Assume that $(p|q):=\dim \DS_{\ol{x}}\DS_y (N)$. Then
$ \dim \DS_{x+y}(N)=(p-j|q-j)$
for some $j\geq 0$. Moreover, $\dim \DS_{x+y}=(p|q)$ if $xN\cap yN=0$.
\end{lem}
\begin{proof}
By Lemma 3.1 of \cite{Sh}, $\DS_{x+y}N\cong DS_{x+y}N^{x^2,y^2}$, so we may assume that $x^2=y^2=0$.  This reduces the statement to the case when $\fg$ is $(0|2)$-dimensional commutative Lie superalgebra and $N$ is a $\mathbb{C}$-graded $\fg$ module
where $x,y$ have different degrees. In particular, for $v\in N$ 
the equality $xv=yv$ implies $xv=0$.

The indecomposable finite-dimensional modules
 over the ring $\mathbb{F}[u,v]/(u^2,v^2)$ were classified in~\cite{Ringel}.
From this classification it follows the indecomposable summands of $N$ are, up to a parity change,  from the following list:
a $4$-dimensional projective modules $M_4$ satisfying
$\DS_x(M_4)=\DS_{x+y}(M_4)=0$ and the ``zigzag'' modules $V_s^{\pm}$.
Each zigzag modules has a basis $\{v_i\}_{i=1}^s$ with $p(v_{i+1})=\ol{i}$; 
we depict each module by the diagram, where $xv_i=v_{i+1}$ is depicted as
$v_i{\longleftarrow}v_{i+1}$ and  $yv_i=v_{i-1}$ is depicted as
 $v_i{\longrightarrow}v_{i-1}$.
 We have
$$\begin{array}{ll}
V_{2n}^+: & v_1{\longrightarrow}v_2\longleftarrow v_3\longrightarrow v_4\longleftarrow\ldots \longrightarrow v_{2n}\\
V_{2n}^-: & v_1{\longleftarrow}v_2\longrightarrow v_3\longleftarrow v_4\longrightarrow\ldots \longleftarrow v_{2n}
\\
V_{2n-1}^+: & v_1{\longrightarrow}v_2\longleftarrow v_3\longrightarrow v_4\longleftarrow\ldots \longleftarrow v_{2n-1}\\
V_{2n-1}^-: & v_1{\longleftarrow}v_2\longrightarrow v_3\longleftarrow v_4\longrightarrow\ldots \longrightarrow v_{2n-1}
\end{array}$$
(The modules $V_1^+\cong V_1^-$ are trivial). One sees that
$$\begin{array}{l}
\dim\DS_{\ol{x}}\DS_y (V^{\pm}_{2n-1})=\dim\DS_y(V^{\pm}_{2n+1})=\dim\DS_{x+y}(V^{\pm}_{2n+1})=(1|0),\\
\DS_y(V^+_{2n})=\DS_y(V^-_2)=0,\ \ \  \DS_{x+y}(V^{\pm}_{2n})=0\end{array}$$
and $\dim \DS_{\ol{x}}\DS_y(V^-_{2n}))=\dim\DS_y(V^-_{2n})=(1|1)$ for $n>1$.
This gives $\dim \DS_{x+y}(N)=(p-j|q-j)$ for some $j\geq 0$. 
If $xN\cap yN=0$, then  $N$ is a direct sum 
and the modules of the form $V_1^{\pm}$, $V^{\pm}_2$, $V^-_3$.
By above, this gives $\dim \DS_{x+y}(N)=\DS_{\ol{x}}\DS_y (N)$.
\end{proof}


\begin{thebibliography}{MMM}


\bibitem{Br} J.~Brundan, {\em Kazhdan-Lusztig polynomials and character formulae
for Lie superalgebra $\fq(n)$}, Adv. Math. {\bf 182} (2004), 28--77.

\bibitem{BrNick} J.~Brundan, N.~Davidson, {\em Type C blocks
of super category $\CO$}, Math. Z. {\bf 293} (2019), no. 3-4, 867–-901.

\bibitem{BS2}  J.~Brundan, C.~Stroppel, 
{\em Highest weight categories arising from Khovanov's diagram algebra. II:  Koszulity}, Transform. Groups, {\bf 15}, (1), (2010).


\bibitem{BS4} J.~Brundan, C.~Stroppel, 
{\em Highest weight categories arising from Khovanov's diagram algebra. IV:    the general linear supergroup}, J. Eur. Math. Soc. (JEMS),
{\bf 14}, (2) (2012).


\bibitem{Cheng} S.~J.~Cheng, {\em Supercharacters of queer Lie superalgebras},
J. Math. Phys. {\bf 58} (2017), no. 6, 061701.


\bibitem{ChK} J.-S.~Cheng, J.-H.~Kwon, {\em Finite-dimensional half-integer weight modules over queer Lie superalgebras} Comm. Math. Phys., {\bf 346} (2016), 945--965.

\bibitem{DS} M.~Duflo, V.~Serganova,  {\em On associated variety for
  Lie superalgebras}, arXiv:0507198.

\bibitem{EAS} I. Entova-Aizenbud, V. Serganova,
{\em Duflo-Serganova functor and superdimension formula for the periplectic Lie superalgebra},
arXiv: {1910.02294},
(2019)

\bibitem{Frisk} A.~Frisk, {\em Typical blocks of the category $\CO$ for the queer
Lie superalgebra}, J. Algebra and Applications, {\bf 6}, No. 5 (2006).


\bibitem{Gq} M.~Gorelik,  {\em Shapovalov determinants for $Q$-type Lie superalgebras}, IMRP (2006), Art. Id. 96895, 1--71.

\bibitem{Gcore} M.~Gorelik, {\em Depths and cores in the light of DS-functors},
arXiv:2010.05721.


\bibitem{Gdex} M.~Gorelik,  {\em Bipartite extension graphs and the Dulfo--Serganova functor}, arXiv: 2010.12817.

\bibitem{Gdexnew} M.~Gorelik,  {\em On modified extension graphs of a fixed atypicality}, 	arXiv:2204.02759.

\bibitem{GH} M.~Gorelik, T.~Heidersdorf, {\em Semisimplicity of the $\DS$ functor for the orthosymplectic Lie superalgebra}, arXiv: 2010.14975.


\bibitem{GHSS} M.~Gorelik, C.~Hoyt, V.~Serganova, A.~Sherman, {\em The Duflo-Serganova functor, vingt ans apr\`es}, arXiv: 2203.00529.

\bibitem{GSS} M.~Gorelik, V.~Serganova, A.~Sherman, {\em Grothendieck rings of quasireductive Lie superalgebras}.

\bibitem{Nicki} N.~Grantcharov, V.~Serganova {\em Extension quiver for Lie superalgebra
$\fq(3)$}, SIGMA Symmetry Integrability Geom. Methods Appl. 16 (2020), Paper No. 141.



\bibitem{GS} C. Gruson, V. Serganova, {\em Cohomology
  of generalized supergrassmanians and character formulae for basic classical Lie superalgebras},
Proc. London Math. Soc., (3), {\bf 101} (2010), 852--892.




\bibitem{GSBGG} C. Gruson, V. Serganova, {\em
  Bernstein-Gelfand-Gelfand reciprocity and
  indecomposable projective modules for classical algebraic
  supergroups}, Mosc. Math. J., {\bf 13} (2013), no. 2, 281–-313.

\bibitem{HW} T.~Heidersdorf, R.~Weissauer {\em Cohomological tensor functors on representations of the
    general linear supergroup}, arXiv:1406.0321, to appear in Mem. Am. Math. Soc.
%

\bibitem{MS} I. Musson, V. Serganova {\em Combinatorics of
  character formulas for the
  Lie superalgebra $\fgl(m|n)$},
Transform. Groups {\bf 16} (2011), no. 2, 555–-578.



\bibitem{P} I.~Penkov, {\em Characters of typical irreducible finite-dimensional
$\fq(n)$-modules}, Funct. Anal. Appl. {\bf 20} (1986), 30--37.
%

\bibitem{PS1} I.~Penkov, V.~Serganova, {\em Characters of finite-dimensional irreducible $\fq(n)$-modules}, Lett.
Math. Phys. {\bf 40 }(1997), 147–-158.
%
\bibitem{PS2} I.~Penkov, V.~Serganova, {\em Characters of irreducible $G$-modules and cohomology of $G/P$ for
the supergroup $G = Q(N)$}, J. Math. Sci., {\bf 84}, (1997), 1382–-1412.

\bibitem{Ringel} C.~M.~Ringel {\em The indecomposable representations of the dihedral 
$2$-groups}, Mathematische Annalen, {\bf 214}, 
(1975), 19--34.

\bibitem{Skw} V.~Serganova {\em On the superdimension of an
  irreducible representation of a
  basic classical Lie superalgebra},
  in Supersymmetry in mathematics and physics, 253-–273, Lecture Notes
  in Math.,
  2027, Springer, Heidelberg, (2011).



\bibitem{S} V.~Serganova {\em Finite-dimensional representation of algebraic supergroups},  Proceedings of the International Congress of Mathematicians—Seoul 2014. Vol. 1, 603–-632, Kyung Moon Sa, Seoul, 2014.




\bibitem{Serq} A.~Sergeev,  {\em The centre of the enveloping algebra
for  Lie superalgebra $Q(n,\mathbb{C})$}, Lett. Math. Phys. {\bf 7} (1983), no. 3, 177--179.

\bibitem{Sh} A.~Sherman, {\em A few remarks on symmetries of the {Duflo-Serganova} functor}, arXiv:2201.12808.

\end{thebibliography}
\end{document}